\documentclass{amsart}
\usepackage{amssymb,amsmath}
\usepackage{amsfonts}
\usepackage{latexsym}
\usepackage{comment,color, marginnote}
\usepackage{mathtools}
\mathtoolsset{showonlyrefs}

\usepackage[left=2.65cm,right=2.65cm,top=2.7cm,bottom=2.7cm]{geometry}
\numberwithin{equation}{section}

\newtheorem{The}{Theorem}[section]
\newtheorem{Lem}[The]{Lemma}
\newtheorem{Pro}[The]{Proposition}

\newtheorem{Cor}[The]{Corollary}

\theoremstyle{remark}
\newtheorem{Rem}[The]{Remark}

\newcommand{\RR}{\mathbb R}
\newcommand{\ZZ}{\mathbb Z}
\newcommand{\TT}{\mathbb T}
\newcommand{\ff}{\mathbf{f}}
\newcommand{\g}{\mathbf{g}}
\newcommand{\dd}{\mathbf{d}}
\newcommand{\hh}{\mathbf{h}}
\newcommand{\jj}{\mathbf{j}}
\newcommand{\LL}{\mathbf{L}}
\newcommand{\kk}{\mathbf{k}}
\newcommand{\Ccal}{\mathcal{C}}
\newcommand{\Ucal}{\mathcal U}

\newcommand{\Scal}{\mathcal{S}}

\newcommand{\Pcal}{\mathcal{P}}

\newcommand{\cqfd}{\mbox{}\nolinebreak\hfill\rule{2mm}{2mm}\medbreak\par}

\renewcommand{\Theta}{\textnormal{O}}
\newcommand{\deps}{d_\varepsilon}

\def\Cc{\mathbf{C}}

\def\mm{{\bf m}}

\begin{document}

\title[Low-energy dynamics in generic potential fields]{Low-energy dynamics in generic potential fields:\\ Hyperbolic periodic orbits and non-ergodicity}
\author{Alberto Enciso, Manuel Garzón, Daniel Peralta-Salas}

\begin{abstract}
We prove that,  on each low energy level, the natural Hamiltonian system defined by a generic smooth potential on $\mathbb T^2$ exhibits an arbitrarily high number of hyperbolic periodic orbits and a positive-measure set of invariant tori. Hence, quasi-periodic motion and hyperbolic behavior typically coexist in the low-energy dynamics of natural Hamiltonian systems with two degrees of freedom.
\end{abstract}

			
			

	
	\maketitle\section{Introduction}

The study of the generic properties of a Hamiltonian system, and in particular the question of whether such systems are typically ergodic or chaotic, is a central topic in dynamical systems, which traces back to the work of Poincar\'e on celestial mechanics and the formulation of the ergodic hypothesis in statistical mechanics due to Boltzmann and Gibbs~\cite{Birk}.

A satisfactory picture of the typical dynamics of a Hamiltonian system did not emerge until the 1970s. Specifically, Markus and Meyer~\cite{MM1,MM2} used transversality theory to prove that a generic (in the Baire category sense) subset of Hamiltonian systems is not integrable, and showed that an open and dense subset of Hamiltonian systems exhibits a positive measure set of invariant tori as a consequence of KAM theory. In particular, a generic Hamiltonian system on a compact symplectic manifold is neither integrable nor ergodic.  

On the contrary, the case of natural Hamiltonian systems is much less studied, in spite of its central relevance in mathematical physics. We recall that, given a compact Riemannian manifold $(M,g)$ and a function $V\in C^\infty(M)$ the associated \textsl{natural Hamiltonian system} on the cotangent bundle $T^*M$ is defined by a Hamiltonian function of the form
    \begin{equation*}
    H_V(p,q):=\frac12 |p|^2_g+V(q)\,,
    \end{equation*}
where the kinetic term is given by the Riemannian metric $g$ and $V$ is the associated potential function. 

Since the kinetic part is fixed in the setting of natural Hamiltonians, and only the potential (a function on $M$) is perturbed, it is difficult to study generic properties of the Hamiltonians $H_V$. The difficulty is that transversality methods fail dramatically when one can only perturb the Hamiltonian by changing a potential function. Markus and Meyer showed~\cite[Theorem~2']{MM2} that for generic potentials the dynamics is not ergodic, and they raised the question of whether the dynamics is integrable or not for a generic potential. If we study the Hamiltonian dynamics for high energy values, i.e., values of $H_V$ close to the maximum value attained by $V$ on $M$, one can effectively use variational methods to prove that, under suitable conditions, the topological entropy of the system is positive~\cite{Ko,BR,BR2}. However, the study of the generic properties of the system for low energy levels is much more complicated~\cite[Section~5]{MM2}. A probabilistic approach to the study of chaos and non-ergodicity for natural Hamiltonians was developed in~\cite{EPR2}, using Gaussian random potentials, but the results do not hold for low energy levels of the Hamiltonian (i.e., those close to the minimum of the potential function). It is also known~\cite{Bernard} that, for a generic potential, there exists {\sl one} low energy level where all periodic orbits are non-degenerate, although no lower bounds are provided for the number of such orbits. Other than that, our understanding of the low-energy particle dynamics in a generic potential field remains rather limited.

For simplicity, in this article we shall only consider the case that $M$ is a $2$-dimensional torus $ \TT^2:= (\RR/\ZZ)^2$  endowed with the flat metric. We can then use global coordinates $(p,q)\in \mathbb R^2\times \mathbb T^2$ to parametrize the phase space of the Hamiltonian, which takes the form
\begin{equation}\label{E.natHam}
    H_V(p,q)=\frac12|p|^2+V(q)\,,
\end{equation}
where $V\in C^\infty(\TT^2)$. The corresponding dynamical system $X_V$ describing the motion of a particle in the potential field $V$ is
\begin{equation}\label{E.HamODE}
    \dot q = p\,,\qquad \dot p = -\nabla V(q) \qquad \Longleftrightarrow \qquad X_V=p\,\partial_q-\nabla V(q)\,\partial_p\,.
\end{equation}
Of course, the function $H_V$ is a first integral (i.e., a conserved quantity) of the dynamics. Since the addition of a constant to the potential $V$ does not change the Hamiltonian dynamics, we can safely normalize~$V$ so that $$\min_{q\in\TT^2} V(q)=0\,.$$

We are interested in the properties of the Hamiltonian system $H_V$ for a generic potential $V$, so we endow the space $C^\infty(\TT^2)$ with the Whitney $C^\infty$-topology. It is well known that this makes $C^\infty(\TT^2)$  a Baire space, meaning that every countable intersection of open and dense sets is also dense. We say that a set of potentials is \textsl{generic} if it is a residual subset of $C^\infty(\TT^2)$ (in the Baire sense).

Our main result shows that, for a generic potential $V$, the Hamiltonian system $H_V$ is non-ergodic and exhibits an arbitrarily large number of hyperbolic periodic orbits on any small enough energy level $h$.

	\begin{The}\label{MainTheorem}
	Let $N_0$ be a positive integer. There exists a generic set $\mathcal V\subset C^\infty(\TT^2)$ and a small constant $\delta_0$ such that the Hamiltonian system $H_V$ with $V\in\mathcal V$ exhibits a positive measure set of invariant tori and at least $N_0$ hyperbolic periodic orbits on each energy level $H^{-1}_V(h)$ for $h\in(0,\delta_0)$.
	\end{The}

To put this result in perspective, let us recall that hyperbolic periodic orbits are the natural starting point to establish the presence of chaos in dynamical systems. On one hand, a theorem of Katok~\cite{Katok} guarantees (in the setting of $3$-dimensional flows, as is the case for the restriction of our Hamiltonian flows to each invariant energy level) that the existence of
hyperbolic periodic orbits is a necessary condition for having positive topological entropy. On the other hand, it is well known that if the stable and unstable manifolds of a hyperbolic orbit intersect transversely, the dynamics must exhibit Smale horseshoes and positive entropy. However, our analysis does {not} enable us to conclude that, for a generic potential, these manifolds intersect transversely: this requires delicate control of the global structure of the invariant manifolds, which we cannot extract from local normal form computations alone because, on low energy levels, it is an exponentially small phenomenon.

Let us briefly discuss the idea of the proof.    
The starting point is the fact that, near a nondegenerate global minimum of the potential, the system is conjugate, via a Birkhoff normal form, to two weakly coupled oscillators. The invariant tori of the linearized system (that is, the decoupled oscillators) play a major role in our analysis of the low-energy dynamics. The non-ergodicity part of the theorem, which builds on the argument sketched by Markus and Meyer in~\cite{MM1, MM2}, is relatively simple and requires little information about the potential. In contrast, proving the existence of {hyperbolic} periodic orbits is considerably more involved and requires much more information on the structure of the normal form. 

The key step of the proof is to compute the coefficients of the normal form very explicitly in terms of derivatives of $V$ at its minimum. In the aforementioned setting, these coefficients and the frequencies of the linearized system determine the twist coefficient explicitly in terms of derivatives of $V$ at the minimum. Since this twist is nonzero under an open and dense nonresonance condition, we can use Moser’s twist theorem to conclude the existence of a positive-measure set of invariant tori on every sufficiently small energy level, which ensures that the Hamiltonian flow is not ergodic.

To obtain the hyperbolic periodic orbits, we analyze resonant tori using a subharmonic Melnikov method. To this end, we pick a dense (but not open) set of potentials with rational frequencies (with quotient, say, $n/m$) and write the $m$-th iterate of the Poincar\'e map as a small  perturbation of a rigid rotation with very {explicit} nonzero resonant coefficients. A careful analysis of this iterate produces $2m$ period-$m$ orbits for the flow, half elliptic and half hyperbolic. With some further elaboration, we conclude that, for a generic set of potentials, on each low energy level there is a fixed but arbitrarily large number of hyperbolic periodic orbits.

The paper is organized as follows. 
In Section~\ref{Sec2} we compute the Birkhoff normal form of the Hamiltonian near the global minimum of the potential, obtaining explicit formulas for the coefficients in terms of the derivatives of $V$. Section~\ref{Sec3-Poincare} introduces a rescaled system to study low-energy dynamics and derives the Poincar\'e map, which is written in terms of a ``reduced'' non-autonomous system of equations in two variables. Section~\ref{Sec4-Nonerg} applies Moser’s twist theorem to establish non-ergodicity of generic potentials at low energy levels. In Section~\ref{S:hyper} we analyze resonant tori to prove the existence of an arbitrarily large number of low-energy elliptic and hyperbolic periodic orbits. The proofs of a number of necessary technical results are relegated to Sections~\ref{Section-ReversedAction}--\ref{Section-Proof2}. Finally, the very short Appendix~\ref{Appendix-Denseclass} shows that the set of potentials satisfying the arithmetic condition needed in Section~\ref{S:hyper} is dense.

	\section{Computation of the Birkhoff normal form at the minimum}\label{Sec2}
	
    Let $\mathcal{V}_0\subset C^\infty(\TT^2)$ be the set of Morse functions {whose critical points take different values}, which is well known to be open and dense. In particular, the unique global minimum of $V\in\mathcal V_0$ (whose value is assumed to be zero) is non-degenerate. In what follows, let us choose the coordinates $q\in\TT^2$ so that the global minimum of $V$ is attained at $q=0$.
    
    For any $V\in\mathcal{V}_0$, we denote by $\{\omega_1^2,\omega_2^2\}\subset(0,\infty)$ the eigenvalues of the Hessian matrix $D^2V(0)$. Then, the eigenvalues of the linearized Hamiltonian vector field $DX_V(0)$ are given by $\{\pm i\omega_1, \pm i\omega_2\}$, with $\omega_k>0$. Moreover, basic linear algebra computations show the existence of a {(unique if $\omega_1\neq \omega_2$)} canonical transformation
	\begin{equation}\label{Canonical}
		{\begin{pmatrix}
			\tilde{q}\\\tilde{p}
		\end{pmatrix}= \begin{pmatrix}
		\mathbf{A}&0\\0&(\mathbf{A}^{-1})^{\intercal}
	\end{pmatrix}\begin{pmatrix}
	q\\p
\end{pmatrix},}
	\end{equation}with $\mathbf{A}$ an isometry, such that, in a {neighborhood of $(0,0)$}, the Hamiltonian $H_V(p,q)$ reads in these variables as
	\begin{equation*}
		\tilde H_V(\tilde p,\tilde q)=\dfrac{|\tilde p|^2}{2}+ \dfrac{\omega^2_1}{2}\tilde q_1^2+\dfrac{\omega^2_2}{2}\tilde q_2^2+\sum\limits_{|\gamma|=3}^N\dfrac{\partial^\gamma \tilde V(0)}{\gamma!}\tilde q^\gamma+H_{\mathrm{R}}(\tilde q).
	\end{equation*}
	{Here $\tilde H_V(\tilde p,\tilde q)=H_V(p,q)$ and $\tilde V(\tilde{q})=V(q)$ are the expressions of the Hamiltonian and of the potential in the new variables, and we are using the standard multi-index notation for $\gamma\in\mathbb{N}^2$.} Also, $H_{\mathrm{R}}$ denotes the remainder term of the $N^{\text{th}}$-order Taylor expansion of $V$ at zero, with $N\in\mathbb{N}$ as big as desired. 
    
    {Since there will be many changes of variable, throughout the paper we will abuse the notation and simply denote the expression of $H_V(p,q)$ and~$V(q)$ in the coordinates $(\tilde p,\tilde q)$ by $H_V(\tilde q,\tilde p)$, $V(\tilde q)$. Furthermore, for the ease of notation we will drop the tildes completely,} so that our Hamiltonian simply reads as
    \begin{equation}\label{Ham}
		 H_V(p,q)=\dfrac{| p|^2}{2}+ \dfrac{\omega^2_1}{2} q_1^2+\dfrac{\omega^2_2}{2} q_2^2+\sum\limits_{|\gamma|=3}^N\dfrac{\partial^\gamma  V(0)}{\gamma!} q^\gamma+H_{\mathrm{R}}( q).
	\end{equation}
    In our coordinates, it thus holds that $\omega_k=[{\partial_{q_k}^2V(0)}]^{1/2}$ for $k\in {1,2}$.
	It is easy to see that for sufficiently small energy values $h>0$, the level set $H_V^{-1}(h)$ is diffeomorphic to the {$3$-sphere (in fact, in our coordinates it is a slightly deformed ellipsoid)}. 	
	
	The following proposition is the main result of this section. It provides the \textsl{Birkhoff normal form} of the Hamiltonian, which will be used to describe the dynamics on low energy levels. {Note that, following the same philosophy as above, we will simply denote by $H_V(I,\varphi)$ the expression of~$H_V(p,q)$ in the variables~$(I,\varphi)$.}
	
	\begin{Pro}\label{Pro-BNF} Let $V\in\mathcal{V}_0$ be such that its eigenvalues satisfy the following non-resonant condition:
	\begin{equation}\label{nonresonant}
		m\omega_1+n\omega_2\neq0,\quad \hbox{for all }m,n\in\mathbb{Z}\setminus\{0\},\ \hbox{with }|m|+|n|\leq4.
	\end{equation}
	Let $\{\beta_1,\beta_2,\beta_3\}$ be the constants
	\begin{align}
		\beta_1&=\frac{\partial_{q_1}^4 V(0)}{8\omega_1^2}+\frac{\left(\partial_{q_1}^3 V(0)\right)^2}{6\omega_1^4}+\frac{\left(\partial_{q_1}^2\partial_{q_2} V(0)\right)^2}{16\omega_1^2}\left[\frac{1}{(2\omega_1+\omega_2)^2}+\frac{1}{\omega_2(\omega_2-2\omega_1)}-\frac{4}{\omega_2^2}\right],\\[1mm]
		\beta_2&=\frac{\partial_{q_2}^4 V(0)}{8\omega_2^2}+\frac{\left(\partial_{q_2}^3 V(0)\right)^2}{6\omega_2^4}+\frac{\left(\partial_{q_1}\partial_{q_2}^2 V(0)\right)^2}{16\omega_2^2}\left[\frac{1}{(2\omega_2+\omega_1)^2}+\frac{1}{\omega_1(\omega_1-2\omega_2)}-\frac{4}{\omega_1^2}\right],\\[1mm]
\beta_{3}&=\frac{\partial_{q_1}^2\partial_{q_2}^2 V(0)}{4\omega_1\omega_2}+\frac{\left(\partial_{q_1}^2\partial_{q_2} V(0)\right)^2}{32\omega_1^2\omega_2}\left[\frac{1}{2\omega_1-\omega_2}+\frac{\omega_1}{2\omega_1+\omega_2}\right]+\frac{\left(\partial_{q_1}\partial_{q_2}^2 V(0)\right)^2}{32\omega_1\omega_2^2}\left[\frac{1}{2\omega_2-\omega_1}+\frac{\omega_2}{2\omega_2+\omega_1}\right].
	\end{align} Then, there exist a set $\mathcal{U}_I=[0,d_1)\times[0,d_2)$, {a neighborhood of the origin }$\mathcal{U}_0\subset\TT^2\times\RR^2$, and a canonical transformation
    \begin{equation*}
        \Phi:\mathcal{U}_I \times\TT^2\rightarrow\mathcal{U}_0,\quad\Phi(I,\varphi)=(p,q),
    \end{equation*}such that the Hamiltonian function \eqref{Ham} reads as:
		\begin{equation}\label{BirkNF}
			H_V\circ \Phi \equiv H_V(I,\varphi)=\sum\limits_{k}\omega_kI_k+\frac{1}{2}\sum\limits_{k}\beta_k I_k^2+\beta_{3}I_1 I_2 +\sum\limits_{l=5}^NH_l(\xi,\eta)+H_{\mathrm{R}}(\xi,\eta),
		\end{equation}
		for any arbitrary natural number $N\geq5$. Moreover, each function $H_l$ is a homogeneous polynomial of degree $l\geq5$ in the variables
        \begin{equation*}
            \xi_k=\sqrt{I_k}e^{-i\varphi_k},\quad\eta_k=i\bar{\xi}_k,\quad\hbox{for }k\in\{1,2\},
        \end{equation*}
        and the remainder satisfies $$H_{\mathrm{R}}=\Theta(I_1^{N/2}+I_2^{N/2}).$$
	\end{Pro}  
	\proof	Let us consider $(q,p)\in \mathcal U_0$. The eigenvectors of $DX_V(0)$ can be parameterized as:
	\begin{equation*}
		\mathbf{v}_{\lambda}=\left\{\begin{array}{lr}
			\left(u,0,\lambda u,0\right)^{\intercal},&\hbox{for any } u\in\mathbb{C},\  \hbox{if }\lambda\in\{\pm i\omega_1\},\\
			\left(0,u,0,\lambda u\right)^{\intercal},&\hbox{for any } u\in\mathbb{C},\  \hbox{if }\lambda\in\{\pm i\omega_2\}.
		\end{array}\right.
	\end{equation*}
	We take the following basis of eigenvectors 
    \begin{equation}
        \mathbf{v}_{-i\omega_1}=\dfrac{1}{\sqrt2}\left(\dfrac{1}{\sqrt{\omega_1}},0,-i\sqrt{\omega_1},0\right)^{\intercal},\quad\hbox{ and }\quad\mathbf{v}_{-i\omega_2}=\dfrac{1}{\sqrt2}\left(0,\dfrac{1}{\sqrt{\omega_2}}0,-i\sqrt{\omega_2}\right)^{\intercal},
    \end{equation}
    and since the vector $-i\bar{\mathbf{v}}_\lambda$ is an eigenvector associated with $\bar{\lambda}$, the other two vectors of the basis are
     \begin{equation}
        \mathbf{v}_{i\omega_1}=-i \overline{\mathbf{v}}_{-i\omega_1},\quad\hbox{ and }\quad\mathbf{v}_{i\omega_2}=-i \overline{\mathbf{v}}_{-i\omega_2}.
    \end{equation}
    Then, we define the symplectic matrix
	\begin{align}\label{eigMatrix}
		\Cc:=(\mathbf{v}_{-i\omega_1}|\mathbf{v}_{-i\omega_2}|\mathbf{v}_{i\omega_1}|\mathbf{v}_{i\omega_2})=\tfrac{1}{\sqrt{2}}\begin{pmatrix}
			\tfrac{1}{\sqrt{\omega_1}}&0&\tfrac{-i}{\sqrt{\omega_1}}&0\\
			0&\tfrac{1}{\sqrt{\omega_2}}&0&\tfrac{-i}{\sqrt{\omega_2}}\\
			-i\sqrt{\omega_1}&0&\sqrt{\omega_1}&0\\
			0&-i\sqrt{\omega_2}&0&\sqrt{\omega_2}
		\end{pmatrix}.
	\end{align}
    Therefore, the transformation $(q,p)^{\intercal}=\Cc(x,y)^{\intercal}$ is canonical, where $(x,y)\in\mathbb{C}^4$ are complex vectors that satisfy the relation $(y_1,y_2)=i(\bar x_1,\bar x_2)$. In these coordinates, the Hamiltonian \eqref{Ham} reads as
	\begin{align}\label{Ham2}
	H_V(x,y)=-i\sum\limits_{k}\omega_kx_ky_k+\sum\limits_{l=3}^N\tilde{H}_l(x,y)+H_{\mathrm{R}}(x,y),
	\end{align}
    where, with some abuse of notation, we are denoting  by $H_V(x,y)$, $H_{\mathrm{R}}(x,y)$, the expression of $H_V$ and $  H_{\mathrm{R}}$ in the new coordinates,
	and where every function $\tilde{H}_l$ is a homogeneous polynomial of degree $l$. In particular, it can be written as:
	\begin{equation}\label{HomPol}
    \tilde{H}_l(x,y)=\sum_{|\gamma|=l}H_\gamma(x,y)=\sum_{(\gamma^1,\gamma^2)\in\Gamma_l}C[\gamma^1,\gamma^2](x_1,y_1)^{\gamma^1}(x_2,y_2)^{\gamma^2},
	\end{equation}
	with
    \begin{equation}\label{set-Gamma_l}
        \Gamma_l=\{(\gamma^1,\gamma^2)\in\mathbb{N}^4:\ |\gamma^1|+|\gamma^2|=l\}.
    \end{equation}Here and in what follows we use the standard multiindex notation. In order to compute $C[\gamma^1,\gamma^2]$, for any vector $\alpha\in\mathbb{N}^2$, we denote by $C(\alpha)$ the binomial coefficient associated with the term $a^{\alpha_1}b^{\alpha_2}$ in the expansion of $(a-ib)^{|\alpha|}$, that is
	\begin{equation}\label{Binomial}
		C(\alpha)=(-i)^{\alpha_2}\dfrac{|\alpha|!}{\alpha!},\quad\hbox{for all } \alpha\in\mathbb{N}^2.
	\end{equation}
	As usual, $\alpha!:=\alpha_1!\alpha_2!$ and $|\alpha|:=\alpha_1+\alpha_2$. Since
	\begin{equation*}
		q_k=\dfrac{1}{\sqrt{2\omega_k}}\left(x_k-iy_k\right),\ \hbox{for }k=1,2,
	\end{equation*}
	it is straightforward to verify that
	\begin{equation}\label{Constant1}
		C[\gamma^1,\gamma^2]=\dfrac{\partial_{q_1}^{|\gamma^1|}\partial_{q_2}^{|\gamma^2|}V(0)}{|\gamma^1|!|\gamma^2|!}\prod\limits_{k=1}^{2}\dfrac{C(\gamma^k)}{\sqrt{(2\omega_k)^{|\gamma^k|}}},\quad \hbox{for all }(\gamma^1,\gamma^2)\in\Gamma_l,\hbox{ with }l\geq3,
	\end{equation}
	which is clearly nonzero for a generic potential $V\in\mathcal{V}_0$.
	
    Now we proceed to summarize the classical method in order to construct canonical transformations between two neighborhoods of the origin in $\mathbb{C}^{4}$, denoted by $\Ucal_1$ and $\Ucal_2$. The reader may check the details e.g. in \cite[Sect. 3 \& 30]{SM}. To avoid repetitions in the proof, any open set $\Ucal$ denotes a neighborhood of the origin in the corresponding space, which is always assumed to be sufficiently small. To this aim, let us consider a function $\nu:\Ucal_\nu\subset\mathbb{C}^4\rightarrow\mathbb{C}$ of the form \begin{equation}\label{generatrix}
		\nu(x,\eta)=x\eta+\sum\limits_{l=3}^{N_0}\nu_l(x,\eta),
	\end{equation}
	where every $\nu_l(x,\eta)$ is again a homogeneous polynomial of degree $l$, and $N_0\in\mathbb{N}$ is arbitrary. Let us also introduce the variables
	\begin{equation}
		\left\{\begin{array}{l}
			\xi_k=\partial_{\eta_k}\nu(x,\eta)=x_k+\sum\limits_{l=3}^{N_0}\partial_{\eta_k}\nu_l(x,\eta),\\
			y_k=\partial_{x_k}\nu(x,\eta)=\eta_k+\sum\limits_{l=3}^{N_0}\partial_{x_k}\nu_l(x,\eta),\qquad k=1,2.
		\end{array}\right.
	\end{equation}
	Since the determinant $|\partial^2_{x_k\eta_j}\nu(x,\eta)|$ is nonzero in $\Ucal_\nu$, the Implicit Function Theorem can be applied with respect to the variable $x$ to the function
    \begin{equation}
        F(x,\xi,\eta)=\xi+\partial_\eta \nu(x,\eta),
    \end{equation}so that there exists an open set $\Ucal_x\subset\mathbb{C}^2$ and a map ${\bf x}:\Ucal_1\to\Ucal_x$ such that
	\begin{equation}
		\left\{\begin{array}{l}
			x_k=\xi_k-\sum\limits_{l=3}^{N_0}\partial_{\eta_k}\nu_l({\bf x}(\xi,\eta),\eta),\\
			y_k=\eta_k+\sum\limits_{l=3}^{N_0}\partial_{x_k}\nu_l({\bf x}(\xi,\eta),\eta),\qquad k=1,2.
		\end{array}\right.
	\end{equation} 
	By the classical theory of normal forms and Hamiltonian systems, the above transformation is canonical for any function $\nu(x,\eta)$ defined as in \eqref{generatrix}. Nevertheless, it is convenient to express these formulas in terms of $(\xi,\eta)$, which leads to the map $(\phi,\psi)=(\phi_1,\phi_2,\psi_1,\psi_2):\Ucal_1\to\Ucal_2$, where
	\begin{equation}\label{Can-Pol}
    x_k=\phi_k(\xi,\eta):=\sum\limits_{l=1}^{N_0}\phi_{k,l}(\xi,\eta)+R_{x_k},\quad
    y_k=\psi_k(\xi,\eta):=\sum\limits_{l=1}^{N_0}\psi_{k,l}(\xi,\eta)+R_{y_k},\quad k=1,2.
	\end{equation}
	We are identifying $(\phi_{k,1},\psi_{k,1})=(\xi_k,\eta_k)$, and $(\phi_{k,2},\psi_{k,2})=(-\partial_{\eta_k}\nu_3,\partial_{x_k}\nu_3)$. Moreover, when $l\geq3$, the functions $(\phi_{k,l},\psi_{k,l})$ are homogeneous polynomials of degree $l$, and they can be written as follows:
	\begin{equation}\label{Coeff1}
		\left\{\begin{array}{l}
			\phi_{k,l}(\xi,\eta)=-\partial_{\eta_k}\nu_{l+1}(\xi,\eta)+\sum\limits_{(\gamma,\tau)\in\Gamma^l}C_{\phi_k}^{\gamma,\tau}(\nabla_\eta^{l}\nu(\xi,\eta))^\gamma(\xi,\eta)^\tau, \\
			\psi_{k,l}(\xi,\eta)=\ \ \partial_{x_k}\nu_{l+1}(\xi,\eta)+\sum\limits_{(\gamma,\tau)\in\Gamma^l}C_{\psi_k}^{\gamma,\tau}(\nabla_\eta^{l}\nu(\xi,\eta))^\gamma(\xi,\eta)^\tau,
		\end{array}\right.
	\end{equation}
	where the vector $\nabla_\eta^l\nu:=\left(\partial_{\eta_1}\nu_3,...,\partial_{\eta_1}\nu_l,\partial_{\eta_2}\nu_3,...,\partial_{\eta_2}\nu_l\right)$, the set $\Gamma^l$ is defined as follows
	\begin{equation}
		\Gamma^l=\{(\gamma,\tau)\in\mathbb{N}^{2l-4}\times\mathbb{N}^4:|\gamma|\geq1\hbox{ and }\sum\limits_{n=1}^{l-2} (n+1)(\gamma_n+\gamma_{n+l-2}) +|\tau|=l\},
	\end{equation}
    and the constants $C^{\gamma,\tau}_{\phi_k}, C^{\gamma,\tau}_{\psi_k}$ depend on the corresponding multi-indexes. Besides, the remainders $R_{x_k}$, $R_{y_k}$ are functions of order $N_0$ with respect to $(\xi,\eta)$. By substituting \eqref{Can-Pol} into \eqref{Ham2}, the Hamiltonian
    \begin{equation}
        H_V(\xi,\eta):=H_V(\phi(\xi,\eta),\psi(\xi,\eta))
    \end{equation}can be written as a new series of homogeneous polynomials:
	\begin{equation}\label{Ham3}
		H_V(\xi,\eta)=H_2(\xi,\eta)+\sum\limits_{l=3}^{N_0}\tilde{H}_l(\xi,\eta)+\tilde{H}_{R}(\xi,\eta),
	\end{equation}
	where $H_2(\xi,\eta)=-i\sum\limits_{k}\omega_k\xi_k\eta_k$ and $\tilde{H}_{R}(\xi,\eta)$ is clearly different than in \eqref{Ham2}. Moreover, for all $l\geq3$:
	\begin{equation}\label{Ham-Pol}
		H_l(\xi,\eta)=\tilde{H}_l(\xi,\eta)-i\sum_{k}\omega_k\left[\partial_{x_k}\nu_{l}(\xi,\eta)\xi_k-\partial_{\eta_k}\nu_{l}(\xi,\eta)\eta_k\right]+\mathcal{P}_l(\xi,\eta),
	\end{equation}
	where the function $\tilde{H}_l$ is defined in \eqref{HomPol} and $\mathcal{P}_l(\xi,\eta)$ denotes a linear combination of $H_2$ and the functions $H_\gamma$ defined in \eqref{HomPol}, with $3\leq |\gamma|\leq l-1$, all evaluated in the variables from \eqref{Can-Pol}-\eqref{Coeff1} with degree smaller than $l$. Concretely, by comparison of exponents it follows that
	\begin{equation}\label{exponentH2}
		H_2(\phi_{1,l_1},\phi_{2,l_2},\psi_{1,l_3},\psi_{2,l_4})\text{ is in }\mathcal{P}_l(\xi,\eta)\text{ if and only if, }l_1+l_3=l_2+l_4=l.
	\end{equation}
	Similarly,
	\begin{equation}\label{exponents}
		H_\gamma(\phi_{1,l_1},\phi_{2,l_2},\psi_{1,l_3},\psi_{2,l_4})\text{ is in }\mathcal{P}_l(\xi,\eta)\text{ if and only if, }\sum_{k=1}^{4}l_k\gamma_k=l.
	\end{equation}
	Observe that no function in $\mathcal{P}_l$ is evaluated at $\{\partial_{x_k}\nu_{l}, \partial_{\eta_k}\nu_{l}: k=1,2\}$, and so $\mathcal{P}_3=0$. In that case,
	\begin{align}
		H_3(\xi,\eta)=\tilde{H}_3(\xi,\eta)-i\sum_{k}\omega_k\left[\partial_{x_k}\nu_3(\xi,\eta)\xi_k-\partial_{\eta_k}\nu_3(\xi,\eta)\eta_k\right].
	\end{align}
	 After comparing terms, and because of \eqref{nonresonant}, the coefficients of $\nu_3$ can be fixed in order to obtain $H_3=0$. Indeed, each of them is proportional to \eqref{Constant1} for the corresponding $(\gamma_1,\gamma_2)\in\Gamma_3$. In other words, every coefficient at $\nu_3$ is determined by a third order derivative of $V$ at zero. Observe that this fixes $\mathcal{P}_4(\xi,\eta)$, which, according to \eqref{exponents}, is of the form,
	 \begin{equation}\label{P4}
	 	H_2(-\partial_{\eta}\nu_3,\partial_{x}\nu_3)+\sum_{|\gamma|=3}H_\gamma(\phi_{1,l_1},\phi_{2,l_2},\psi_{1,l_3},\psi_{2,l_4}).
	 \end{equation}
	 Notice that at each $H_\gamma$, necessarily it must be $l_j=2$ and $\gamma_j=1$, for some $j$, while the remaining three variables belong to $(\xi,\eta)$, such that at least one of them has null exponent. Furthermore, because of the choice of $\nu_3$, all the constants at $\mathcal{P}_4(\xi,\eta)$ are second order polynomials with respect to the third order derivatives of $V$ at the origin. Finally, it is immediate to see that, for any $\nu_{4}$, the powers and the product between the pairs $\{(\xi_1,\eta_1), (\xi_2,\eta_2)\}$ vanish in the sum 
	 \begin{align}\label{QuadraticTerms}
	 i\sum_{k}\omega_k\left[\partial_{x_k}\nu_{4}(\xi,\eta)\xi_k-\partial_{\eta_k}\nu_{4}(\xi,\eta)\eta_k\right],
	 \end{align} 
	 and so these terms in $H_4$ cannot be canceled under any choice of $\nu_4$. This means that its simplest expression after reduction is
	\begin{equation}
		H_4(\xi,\eta)=-\tfrac{1}{2}\sum\limits_{k}\beta_k\xi_k^2\eta_k^2-\beta_{3}\prod_{k}\xi_k\eta_k,
	\end{equation}
	for which assumption \eqref{nonresonant} is again necessary. Proceeding as above, by \eqref{Ham-Pol}, every $\beta_k$ is determined by the corresponding coefficient of $\tilde{H}_4$ and those from \eqref{P4}. Consequently, each of them is given by a fourth order derivative of $V$ at zero, according to \eqref{Constant1}, plus a quadratic polynomial whose variables are third order derivatives of $V$ at zero. Computing them explicitly, we obtain the expressions defined in the statement. To conclude, \eqref{BirkNF} is obtained after writing $(\xi,\eta)$ in their exponential form, which decomposes $\Ucal_1$ into $\Ucal_I\times\mathbb{T}^2$, which completes the proof.
	\cqfd
	\begin{Rem}
		Observe that the condition \eqref{nonresonant} is open and dense in $ C^\infty(\TT^2)$, hence Proposition~\ref{Pro-BNF} holds for an open and dense set of potentials. Additionally, the integrable part of the Birkhoff normal form \eqref{BirkNF} can be expanded into a series of products and powers of the pairs $\{(\xi_1,\eta_1), (\xi_2,\eta_2)\}$, up to any finite order, as large as desired. The procedure is analogous and only requires extending the non-resonance condition \eqref{nonresonant} to higher orders. More precisely, it follows directly from the fact that, under any choice of the generating function \eqref{generatrix}, the monomials in the pairs $\{(\xi_1,\eta_1), (\xi_2,\eta_2)\}$ always vanish in the expression
		\begin{align}
			i\sum_{k}\omega_k\left[\partial_{x_k}\nu_{l}(\xi,\eta)\xi_k-\partial_{\eta_k}\nu_{l}(\xi,\eta)\eta_k\right]
		\end{align}
		whenever $l\in\mathbb{N}$ is even, and therefore these are precisely the terms that cannot be canceled in \eqref{BirkNF}.
	\end{Rem}

    The following corollary, which is a straightforward consequence of Proposition~\ref{nonresonant}, provides a precise formula for the Birkhoff normal form computed above. In particular, it shows the explicit dependence of the different constants with respect to the potential $V$.
    
    \begin{Cor}\label{Cor-SymplecticTransformation}
    Let $V\in\mathcal{V}_0$ satisfy the hypotheses of Proposition~\ref{Pro-BNF} and let $N\geq5$ be an arbitrary natural number. Then, there exist a set $\mathcal{U}_I=[0,d_1)\times[0,d_2)$, {a neighborhood of the origin }$\mathcal{U}_0\subset\TT^2\times\RR^2$, and a canonical transformation
    \begin{equation*}
        \Phi:\mathcal{U}_I \times\TT^2\rightarrow\mathcal{U}_0,\quad\Phi(I,\varphi)=(p,q),
    \end{equation*}such that, for all $5\leq l\leq N$, the $l^{th}$-order term in the Birkhoff normal form \eqref{BirkNF} is characterized as follows:
    \begin{equation}\label{Hom-Pol-FinalForm}
		H_l(\xi,\eta)=\sum\limits_{\gamma=(\gamma^1,\gamma^2)\in\Gamma_l} C_{\gamma}\ (\xi_1,\eta_1)^{\gamma^1}(\xi_2,\eta_2)^{\gamma^2},
	\end{equation}
	with
		\begin{equation}\label{Constant-FinalForm}
			C_\gamma:=\dfrac{\partial_{q_1}^{|\gamma^1|}\partial_{q_2}^{|\gamma^2|}V(0)}{|\gamma^1|!|\gamma^2|!}\prod\limits_{k=1}^{2}\dfrac{C(\gamma^k)}{\sqrt{(2\omega_k)^{|\gamma^k|}}}+P_\gamma,
		\end{equation}
		where $C(\cdot)$ is the binomial constant from \eqref{Binomial} and $P_\gamma$ is a polynomial of order $\lfloor\tfrac{l}{2}\rfloor$, with respect to the $k^{\mathbf{th}}$-order derivatives of $V$ at zero, for $3\leq k\leq l-1$. Furthermore, for all finite $l\geq5$, the condition
		\begin{equation}\label{NonzeroConstantsGeneric}
			C_\gamma\neq0,\ \hbox{ for all }\gamma\in\Gamma_k,\hbox{ with }5\leq k\leq l,
		\end{equation}
		is open and dense in $ C^\infty(\TT^2)$.
    \end{Cor}
    \proof
    Let the generating function $\nu: \Ucal_\nu\subset\mathbb{C}^4\to\mathbb{C}$ from \eqref{generatrix} be of the form
    \begin{equation}
        \nu(\xi,\eta)=\xi\eta+\nu_3(\xi,\eta)+\nu_4(\xi,\eta),
    \end{equation}
    where the functions $\nu_3$, $\nu_4$ are introduced in the proof of Proposition~\ref{Pro-BNF}, which lead to the Birkhoff normal form \eqref{BirkNF}. In particular, $\nu_l=0$, for every $l\geq5$. In this case, the canonical transformation \eqref{Can-Pol} is of the form
		\begin{equation}\label{canonical-Hyp}
			x_k=\sum\limits_{l=1}^{N_0}\phi_{k,l}(\xi,\eta)+R_{x_k},\ \
			y_k=\sum\limits_{l=1}^{N_0}\psi_{k,l}(\xi,\eta)+R_{y_k},\quad k=1,2,
		\end{equation}
		where the functions $\{\phi_{k,l}$, $\psi_{k,l}\}$ from \eqref{Coeff1} are now characterized as follows. On the one hand, 
		\begin{equation*}
			\left\{\begin{array}{l}
				\phi_{k,3}(\xi,\eta)=-\partial_{\eta_k}\nu_4(\xi,\eta)+\sum\limits_{(\gamma,\tau)\in\Gamma^3}C_{\phi_k}^{\gamma,\tau}(\partial_\eta\nu_3(\xi,\eta))^{\gamma}(\xi,\eta)^\tau, \\
				\psi_{k,3}(\xi,\eta)=\ \ \partial_{x_k}\nu_4(\xi,\eta)+\sum\limits_{(\gamma,\tau)\in\Gamma^3}C_{\psi_k}^{\gamma,\tau}(\partial_\eta\nu_3(\xi,\eta))^{\gamma}(\xi,\eta)^\tau,\quad k=1,2,
			\end{array}\right.
		\end{equation*}
		with $\Gamma^3=\{(\gamma,\tau)\in\mathbb{N}^{2}\times\mathbb{N}^4:|\gamma|=|\tau|=1\}$. Moreover, for all $l\geq4$:
		\begin{equation*}
			\left\{\begin{array}{l}
				\phi_{k,l}(\xi,\eta)=\sum\limits_{(\gamma,\tau)\in\Gamma^l}C_{\phi_k}^{\gamma,\tau}(\nabla_\eta\nu(\xi,\eta))^{\gamma}(\xi,\eta)^\tau, \\
				\psi_{k,l}(\xi,\eta)=\sum\limits_{(\gamma,\tau)\in\Gamma^l}C_{\psi_k}^{\gamma,\tau}(\nabla_\eta\nu(\xi,\eta))^{\gamma}(\xi,\eta)^\tau,\quad k=1,2,
			\end{array}\right.
		\end{equation*}
		with $\nabla_\eta\nu=(\partial_{\eta_1}\nu_3,\partial_{\eta_1}\nu_4,\partial_{\eta_2}\nu_3,\partial_{\eta_2}\nu_4)$, and
		\begin{equation*}
			\Gamma^l=\{(\gamma,\tau)\in\mathbb{N}^{4}\times\mathbb{N}^4:|\gamma|\geq1\hbox{ and }\sum\limits_{n=1}^{2} (n+1)(\gamma_n+\gamma_{n+2}) +|\tau|=l\},
		\end{equation*}
		for $l\geq4$. With this choice, \eqref{Ham-Pol} reduces to $$H_l(\xi,\eta)=\tilde{H}_l(\xi,\eta)+\mathcal{P}_l(\xi,\eta),$$ for all $l\geq5$. In addition, $\mathcal{P}_l(\xi,\eta)$ is a linear combination of $H_2$ and the functions $H_\gamma$ from \eqref{HomPol}, with $3\leq |\gamma|\leq l-1$, all evaluating in the variables \eqref{canonical-Hyp}, according to the corresponding requirements \eqref{exponentH2}-\eqref{exponents}. Then, similarly to \eqref{HomPol}, the homogeneous polynomial $H_l$ can be written as follows
	\begin{equation}
		H_l(\xi,\eta)=\sum\limits_{(\gamma^1,\gamma^2)\in\Gamma_l} C_{\gamma} (\xi_1,\eta_1)^{\gamma^1}(\xi_2,\eta_2)^{\gamma^2},
	\end{equation}
	where at each term, the constant $C_\gamma$ is the sum of the corresponding constant from the polynomials $\tilde{H}_l$ and $\Pcal_l$, i.e.,
		\begin{equation}
        C_\gamma=C[\gamma^1,\gamma^2]+P_\gamma=\dfrac{\partial_{q_1}^{|\gamma^1|}\partial_{q_2}^{|\gamma^2|}V(0)}{|\gamma^1|!|\gamma^2|!}\prod\limits_{k=1}^{2}\dfrac{C(\gamma^k)}{\sqrt{(2\omega_k)^{|\gamma^k|}}}+P_\gamma,
		\end{equation}
        where $C[\gamma_1,\gamma_2]$ is defined in \eqref{Constant1} and $P_\gamma$ is a polynomial of order $\lfloor\tfrac{l}{2}\rfloor$, with respect to the $k^{\mathbf{th}}$-order derivatives of $V$ at zero, for $3\leq k\leq l-1$. \\
        Finally, since the constants $C(\alpha)$ defined in \eqref{Binomial} never vanish, each $C_\gamma$ is an algebraic function on the variables $\{\partial^{\nu}_qV(0):\ 2\leq|\nu|\leq|\gamma|\}$. It then easily follows that the property
		\eqref{NonzeroConstantsGeneric}	is open and dense in $ C^\infty(\TT^2)$, for all finite $l\geq5$.
    \cqfd

\section{Computation of the Poincar\'e map}\label{Sec3-Poincare}
For any $\varepsilon>0$, let us consider the following change of variables in $\Ucal_I$:
\begin{equation}\label{rescaling}
	I=(I_1,I_2)=\varepsilon^2(\tilde{I}_1,\tilde{I}_2)=\varepsilon^2\tilde{I}.
\end{equation} With this substitution, \eqref{BirkNF} takes the form
\begin{equation}
    H_V(\varepsilon^2\tilde{I},\varphi)=\varepsilon^2\sum_{k=1}^2\omega_k \tilde{I}_k+\varepsilon^4H_1(\tilde{I})+\sum_{l\geq5}^{N}\varepsilon^{l}H_{l}(\tilde{I},\varphi)+H_{\mathrm{R}}(\varepsilon^2\tilde{I},\varphi),
\end{equation}
where $H_1$ denotes the quadratic part in \eqref{BirkNF},
\begin{align*}
	H_{l}(\tilde{I},\varphi):=H_{l}(\xi(\tilde{I},\varphi),\eta(\tilde{I},\varphi)),\quad\hbox{for } 5\leq l\leq N,
\end{align*}and
\begin{equation}
	H_{\mathrm{R}}(\varepsilon^2\tilde{I},\varphi):=H_{\mathrm{R}}(\xi(\varepsilon^2\tilde{I},\varphi),\eta(\varepsilon^2\tilde{I},\varphi)).
\end{equation}
In addition, it is convenient to introduce the rescaled Hamiltonian 
\begin{align}\label{BirkNF2}
	\tilde{H}_V(\tilde{I},\varphi):=\varepsilon^{-2}H_V(\varepsilon^2\tilde{I},\varphi)=\sum_{k=1}^2\omega_k \tilde{I}_k+\varepsilon^2H_1(\tilde{I})+\sum_{l\geq5}^{N}\varepsilon^{l-2}H_{l}(\tilde{I},\varphi)+H_{\mathrm{R}}(\tilde{I},\varphi,\varepsilon),
\end{align}
where $H_{\mathrm{R}}(\tilde{I},\varphi,\varepsilon):=\varepsilon^{-2}H_{\mathrm{R}}(\varepsilon^2\tilde{I},\varphi)$, so the remainder is of order $\Theta(\varepsilon^{N-2})$. 
\begin{Rem}\label{Remark-EnergyLevelSets}
    Recall that the Birkhoff normal form \eqref{BirkNF} is valid close to the origin, i.e., on the energy levels
    \begin{equation}
        H_V^{-1}(h),\quad \hbox{with }h\in(0,\delta_0),
    \end{equation}
    for some small enough $\delta_0>0$. On the other hand, for each $\varepsilon>0$, the rescaled Birkhoff normal form \eqref{BirkNF2} holds on the level sets
    \begin{equation}
        \tilde{H}_V^{-1}(h),\quad \hbox{with }h\in(0,h_0),
    \end{equation}
    where $h_0=\varepsilon^{-2}\delta_0$. Nevertheless, we shall see that the constructions developed later in this work require $h_0$ to remain a small ($\varepsilon$-independent) constant. 
\end{Rem}

From now on, we shall work with the rescaled Hamiltonian \eqref{BirkNF2}, dropping the tildes  for notational simplicity. Then, the associated Hamiltonian system  reads
\begin{equation}\label{HamiltonianSystem}
	\left\{\begin{array}{l}
		\dot{I}_k=-\sum\limits_{l=5}^{N}\varepsilon^{l-2}\partial_{\varphi_k}H_{l}(I,\varphi)+\partial_{\varphi_k}H_{\mathrm{R}}(I,\varphi,\varepsilon),\quad k=1,2,\\
		\dot{\varphi}_k=\omega_k+\varepsilon^2\partial_{I_k}H_1(I)+\sum\limits_{l=5}^{N}\varepsilon^{l-2}\partial_{I_k}H_{l}(I,\varphi)+\partial_{I_k}H_{\mathrm{R}}(I,\varphi,\varepsilon).
	\end{array}
	\right.
\end{equation}
When the action variables in Equation~\eqref{HamiltonianSystem} are conserved, the system is completely integrable and every energy level set is densely filled by invariant tori (resonant or nonresonant). We claim that this picture essentially holds for small values of $\varepsilon$, for an open and dense set of potentials $V\in\mathcal{V}_0$. In order to prove this statement, we first show that the dynamics on each energy level can be reduced to a (nonautonomous) planar system. In the following proposition it is convenient to change the notation: $$(I_1,I_2)=(I,J)=\vec{I},\quad (\varphi_1,\varphi_2)=(\varphi,\theta)=\vec{\varphi}.$$
\begin{Pro}\label{Pro-Continuous}
Assume that the potential $V\in\mathcal{V}_0$ satisfies~\eqref{nonresonant}. Then, on every level set of the Hamiltonian~\eqref{BirkNF2} with small enough energy $h>0$, \eqref{HamiltonianSystem} can be reduced to the planar system
\begin{equation}\label{Continuous0}
	\left\{\begin{array}{l}
		\dfrac{dI}{d\theta}=-\dfrac{\varepsilon^3}{\omega_2}\partial_{\varphi}H_5(I,J_0,\vec{\varphi})+\Theta(\varepsilon^4),\\
		\dfrac{d\varphi}{d\theta}=\dfrac{\omega_1}{\omega_2}+\varepsilon^2\Lambda(I)+\dfrac{\varepsilon^3}{\omega_2^2}\left[\omega_2\partial_{I}-\omega_1\partial_{J}\right]H_5(I,J_0,\vec{\varphi})+\Theta(\varepsilon^4),
	\end{array}\right.
\end{equation}
with $J_0:=\tfrac{h-\omega_1 I}{\omega_2}$, and
\begin{equation}
\Lambda(I):=\dfrac{\omega_2(\beta_1I+\beta_{3}J_0)-\omega_1(\beta_2 J_0+\beta_{3}I)}{\omega_2^2}.
\end{equation}
\end{Pro}
\proof
It is clear that for $V$ as in the statement, there exist $\bar{\varepsilon}>0$, and a new set $\mathcal{U}=[0,d_I)\times[0,d_J)$, with $d_I, d_J$ strictly positive and depending on $\bar{\varepsilon}$, such that
\begin{equation}\label{non-degeneracy}
	\partial_{J} H_V(\vec{I},\vec{\varphi})=\omega_2+\varepsilon^2(\beta_2J+\beta_{3}I)+\Theta(\varepsilon^3)>0,
\end{equation}
for all $(\vec{I},\vec{\varphi})\in\mathcal{U}\times\TT^2$, and $\varepsilon<\bar{\varepsilon}$. Moreover, fixing a small enough $h>0$, the level set $$H_V(\vec{I},\vec{\varphi})=h$$ is contained in $\mathcal{U}\times\TT^2$, for every $\varepsilon<\bar{\varepsilon}$. Then, the Implicit Function Theorem gives the existence of a smooth map $J:[0,d_I)\times\TT^2\times[0,\bar{\varepsilon})\rightarrow[0,d_J)$ such that 
\begin{equation}\label{levelSet}
	H_V(I,J(I,\vec{\varphi},\varepsilon;h),\vec{\varphi})=h.
\end{equation}
From now on, we shall write $J:=J(I,\vec{\varphi},\varepsilon;h)$ to simplify the notation. On the other hand, the non-degeneracy condition \eqref{non-degeneracy} also allows us to invert the time variable with respect to $\theta$, transforming \eqref{HamiltonianSystem} on the level set \eqref{levelSet} into the following system:
\begin{equation}\label{Continuous1}
	\left\{\begin{array}{l}
		\vspace{0.5mm}\dfrac{dI}{d\theta}=-\sum\limits_{l=5}^{N}\varepsilon^{l-2}\dfrac{\partial_{\varphi}H_{l}}{\partial_{J}H_V}(I,J,\vec{\varphi},\varepsilon)+\dfrac{\partial_{\varphi}H_{\mathrm{R}}}{\partial_{J}H_V}(I,J,\vec{\varphi},\varepsilon),\\
		\dfrac{d\varphi}{d\theta}=\dfrac{\partial_{I}H_V}{\partial_{J}H_V}(I,J,\vec{\varphi},\varepsilon).
	\end{array}\right.
\end{equation}
In order to expand both equations with respect to $\varepsilon$, we first need to study the function~$J$. To this aim, we compute its Taylor expansion at $\varepsilon=0$ up to order $n<N$:
\begin{align}
	J= \dfrac{h-\omega_1 I}{\omega_2}+\sum\limits_{k=1}^n\dfrac{\varepsilon^k}{k!} \partial^k_\varepsilon J(I,\vec{\varphi},0)+\Theta(\varepsilon^{n+1})=:J_0+\sum\limits_{k=1}^n\varepsilon^k J_k+\Theta(\varepsilon^{n+1}).
\end{align}
Substituting this expression into Equation~\eqref{levelSet}, and using that $\omega_1 I+ \omega_2J_0=h$, we obtain
\begin{equation}\label{levelSet2}
	0=\omega_2\sum\limits_{k=1}^n\varepsilon^k J_k+\varepsilon^2 H_1(I,J)+\sum\limits_{l=5}^{n}\varepsilon^{l-2}H_{l}(I,J,\vec{\varphi})+\Theta(\varepsilon^{n+1}),
\end{equation}
which necessarily implies $J_1\equiv0$. Then, a direct computation gives
$H_1(I,J)=H_1(I,J_0)+\Theta(\varepsilon^2),$ and $H_{l}(I,J,\vec{\varphi})=H_{l}(I,J_0,\vec{\varphi})+\Theta(\varepsilon^2)$, for any $l\geq5,$ which leads to
\begin{equation}
\dfrac{dI}{d\theta}=-\dfrac{\varepsilon^3}{\omega_2}\partial_{\varphi}H_5(I,J_0,\vec{\varphi})+\Theta(\varepsilon^4).
\end{equation}
On the other hand, 
\begin{align}
	\dfrac{d\varphi}{d\theta}=\dfrac{\partial_{I}H_V}{\partial_{J}H_V}(I,J,\vec{\varphi},\varepsilon)=&\ \dfrac{\omega_1+\varepsilon^2\partial_{I}H_1(I,J_0)}{\partial_{J}H_V(I,J,\vec{\varphi},\varepsilon)}+\dfrac{\varepsilon^3}{\omega_2}\partial_{I}H_5(I,J_0,\vec{\varphi})+\Theta(\varepsilon^4)\\
	=:&\ F(\varepsilon)+\dfrac{\varepsilon^3}{\omega_2}\partial_{I}H_5(I,J_0,\vec{\varphi})+\Theta(\varepsilon^4),
\end{align}
so it remains to study the function $F(\varepsilon)$. Let us set $G(\varepsilon):=\partial_{J}H_V(I,J,\vec{\varphi},\varepsilon)$. Then, omitting the arguments $(I,J,\vec{\varphi},\varepsilon)$, it follows that
\begin{align}
	G'(\varepsilon)&=\partial_\varepsilon\partial_{J}H_V+\partial^2_{J}H_V\partial_\varepsilon J,\\
	G''(\varepsilon)&=\partial^2_\varepsilon\partial_{J}H_V+2\partial_\varepsilon\partial^2_{J}H_V\partial_\varepsilon J+\partial^3_{J}H_V(\partial_\varepsilon J)^2+\partial^2_{J}H_V\partial^2_\varepsilon J,
\end{align}and
\begin{align}
G'''(\varepsilon)=&\ \partial^3_\varepsilon\partial_{J}H_V+3\partial^2_\varepsilon\partial^2_{J}H_V\partial_\varepsilon J+3\partial_\varepsilon\partial^3_{J}H_V(\partial_\varepsilon J)^2+3\partial_\varepsilon\partial^2_{J}H_V\partial^2_\varepsilon J\\&+\partial^4_{J}H_V(\partial_\varepsilon J)^3+3\partial^3_{J}H_V\partial_\varepsilon J\partial^2_\varepsilon J+\partial^2_{J}H_V\partial^3_\varepsilon J.
\end{align}
 On the other hand,
\begin{align}
F'(\varepsilon)=&2\dfrac{\varepsilon\partial_{I}H_1}{G}-\dfrac{\omega_1+\varepsilon^2\partial_{I}H_1}{G^2}G',\\
F''(\varepsilon)=&2\dfrac{\partial_{I}H_1}{G}-4\dfrac{\varepsilon\partial_{I}H_1}{G^2}G'-\dfrac{\omega_1+\varepsilon^2\partial_{I}H_1}{G^2}G''+2\dfrac{\omega_1+\varepsilon^2\partial_{I}H_1}{G^3}(G')^2,\\
F'''(\varepsilon)=&-6\dfrac{\partial_{I}H_1}{G^2}G'-6\dfrac{\varepsilon\partial_{I}H_1}{G^2}G''+12\dfrac{\varepsilon\partial_{I}H_1}{G^3}(G')^2-\dfrac{\omega_1+\varepsilon^2\partial_{I}H_1}{G^2}G'''\\&+6\dfrac{\omega_1+\varepsilon^2\partial_{I}H_1}{G^3}G'G''-6\dfrac{\omega_1+\varepsilon^2\partial_{I}H_1}{G^4}(G')^3.
\end{align}
Then, since we have $$G'(0)=0,\quad G''(0)=2\partial_{J}H_1(I,J_0),\quad G'''(0)=3\partial_{J}H_5(I,J_0,\vec{\varphi}),$$ it follows that
\begin{align}
	F'(0)&=0,\quad F''(0)=2\dfrac{\omega_2\partial_{I}H_1(I,J_0)-\omega_1\partial_{J}H_1(I,J_0)}{\omega_2^2},\\
	F'''(0)&=-3\dfrac{\omega_1}{\omega_2^2}\partial_{J}H_5(I,J_0,\vec{\varphi}),
\end{align}
which gives the desired expression \eqref{Continuous0} after substituting into Equation~\eqref{Continuous1}.
\cqfd

\begin{Rem}\label{Twist}
	Notice that $\Lambda(I)$ is a linear function with derivative
	\begin{equation*}
		\Lambda'(I)=-\dfrac{1}{\omega_2^3}\begin{vmatrix}
			\beta_1&\beta_{3}&\omega_1\\ \beta_{3}&\beta_2&\omega_2\\ \omega_1&\omega_2&0
		\end{vmatrix}=:\dfrac{1}{\omega_2^3}\ \Lambda_V,
	\end{equation*}
	which only depends on the derivatives of the potential at zero up to order $4$. The constant $\Lambda_V$ is the so-called \textnormal{twist coefficient} of \eqref{BirkNF} and it is an algebraic function with respect to the derivatives $\{\partial^{\gamma}_q V(0):\ 2\leq|\gamma|\leq4 \}$. Because of that, the condition
	\begin{equation}\label{Isochronic}
		\Lambda_V\neq0,
	\end{equation}
    commonly known as \textnormal{isochronic non-degeneracy} of $V$, is open and dense in $ C^\infty(\TT^2)$.
\end{Rem}
Moreover, since $\partial_I H (\vec{I},\vec{\varphi})$ satisfies a positivity condition analogous to \eqref{non-degeneracy}, the Hamiltonian system \eqref{HamiltonianSystem} can also be reduced to
\begin{equation}
	\left\{\begin{array}{l}
		\dfrac{dJ}{d\varphi}=-\dfrac{\varepsilon^3}{\omega_1}\partial_{\theta}H_5(I_0,J,\vec{\varphi})+\Theta(\varepsilon^4),\\
		\dfrac{d\theta}{d\varphi}=\dfrac{\omega_2}{\omega_1}+\varepsilon^2\Lambda_J(I)+\dfrac{\varepsilon^3}{\omega_1^2}\left[\omega_1\partial_{J}-\omega_2\partial_{I}\right]H_5(L^0,J,\vec{\varphi})+\Theta(\varepsilon^4),
	\end{array}\right.
\end{equation}
with $I_0=\tfrac{h-\omega_2 J}{\omega_1}$, and
\begin{equation}
	\Lambda_J(J)=\dfrac{\omega_1(\beta_2J+\beta_{3}I_0)-\omega_2(\beta_1 I_0+\beta_{3}J)}{\omega_1^2},
\end{equation}
which leads to $\Lambda_J'(J)=\omega_1^{-3}\Lambda_V$. Therefore, it is not relevant which action is used in order to obtain the twist condition~\eqref{Isochronic}.

The following technical lemma concerning the function $\Lambda$ will be key in Section~\ref{S:hyper} in order to study the hyperbolic periodic orbits of the system.

\begin{Lem}
	Let $V\in\mathcal{V}_0$ be a potential which satisfies the non-resonant condition~\eqref{nonresonant} and the twist condition~\eqref{Isochronic}. Then there exists a constant $h^0>0$ such that, for every level set of the Hamiltonian~\eqref{BirkNF2} with energy $h\in(0,h^0)$, the function $\Lambda(\cdot)$ vanishes at some point $I^*_0\in(0,d_I)$.
\end{Lem}
\proof 
First, we notice that $\Lambda(\cdot)$ vanishes at some point $I^*_0\in(0,d_I)$ if and only if,
\begin{equation}\label{ResonanceCondition}
	I^*_0=h\dfrac{\omega_1\beta_2-\omega_2\beta_{3}}{\Lambda_V}\in(0,d_I),
\end{equation}
which holds for any small enough $h>0$ if and only if,
\begin{equation*}
	(\omega_1\beta_2-\omega_2\beta_{3})\Lambda_V^{-1}>0.
\end{equation*} 
Similarly, the function $\Lambda_J(\cdot)$ vanishes at some point $J^*_0\in(0,d_J)$ only when $(\omega_2\beta_1-\omega_1\beta_{3})\Lambda_V^{-1}>0.$ However, since
\begin{equation}
	\omega_1\dfrac{\omega_1\beta_2-\omega_2\beta_{3}}{\Lambda_V}+\omega_2\dfrac{\omega_2\beta_1-\omega_1\beta_{3}}{\Lambda_V}=1,
\end{equation}
at least one of these two conditions holds, thus implying that either $\Lambda$ or $\Lambda_V$ vanishes at some point in their respective domains. Therefore, without any loss of generality, we can assume that condition \eqref{ResonanceCondition} holds, which completes the proof of the lemma.
\cqfd
Finally, let us compute the Poincaré map associated with the reduced system~\eqref{Continuous0}, which is an area-preserving map of the annulus. Specifically, for any  $\theta_0\in\TT$, the $\theta$-integration of \eqref{Continuous0} over $\TT$ yields the map $\mathcal{P}_\varepsilon:(0,d_I)\times\TT\rightarrow(0,d_I)\times\TT$:
\begin{equation}\label{Discrete0}
	\left\{\begin{array}{l}
		I(\theta_0+2\pi)=I_0+f(I,\varphi,\theta_0,\varepsilon),\\
		\varphi(\theta_0+2\pi)=\varphi_0+2\pi\Upsilon(I_0)+g(I,\varphi,\theta_0,\varepsilon),
	\end{array}
	\right.
\end{equation}
where $(I(\theta_0),\varphi(\theta_0))=(I_0,\varphi_0)\in(0,d_I)\times\TT$ is the initial datum. In this expression we have set $\Upsilon(I):=\frac{\omega_1}{\omega_2}+\varepsilon^2\Lambda(I)$, and
\begin{align}
	&f(I,\varphi,\theta_0,\varepsilon):=-\dfrac{\varepsilon^3}{\omega_2}\int_{\theta_0}^{\theta_0+2\pi}\partial_{\varphi_1}H_5(I_0,J_0(I_0),\varphi_0,\theta)d\theta+\Theta(\varepsilon^4),\\
	&g(I,\varphi,\theta_0,\varepsilon):=\dfrac{\varepsilon^3}{\omega_2^2}\int_{\theta_0}^{\theta_0+2\pi}\left[\omega_2\partial_{I}-\omega_1\partial_{J}\right]H_5(I_0,J_0(I_0),\varphi_0,\theta)d\theta+\Theta(\varepsilon^4),
\end{align}
which are well defined provided that $\varepsilon<\bar{\varepsilon}$. 
Similarly, for any $m\in\mathbb{N}$, the $m$-iteration of the map~\eqref{Discrete0} yields the map $\mathcal{P}^m_\varepsilon:(0,d_I)\times\TT\rightarrow(0,d_I)\times\TT$:
\begin{equation}\label{DiscreteM}
	\left\{\begin{array}{l}
		I^m(\theta_0):=I(\theta_0+2\pi m)=I_0+f^m(I,\varphi,\theta_0,\varepsilon),\\
		\varphi^m(\theta_0):=\varphi(\theta_0+2\pi m)=\varphi_0+2\pi m\Upsilon(I_0)+g^m(I,\varphi,\theta_0,\varepsilon),
	\end{array}
	\right.
\end{equation}
with
\begin{align*}
	&f^m(I,\varphi,\theta_0,\varepsilon):=-\dfrac{\varepsilon^3}{\omega_2}\int_{\theta_0}^{\theta_0+m2\pi}\partial_{\varphi_1}H_5(I_0,J_0(I_0),\varphi_0,\theta)d\theta+\Theta(\varepsilon^4),\\
	&g^m(I,\varphi,\theta_0,\varepsilon):=\dfrac{\varepsilon^3}{\omega_2^2}\int_{\theta_0}^{\theta_0+m2\pi}\left[\omega_2\partial_{I}-\omega_1\partial_{J}\right]H_5(I_0,J_0(I_0),\varphi_0,\theta)d\theta+\Theta(\varepsilon^4).
\end{align*}
It is clear that any $2\pi$-periodic solution of \eqref{Continuous0} corresponds to a fixed point of \eqref{Discrete0}, while sub-harmonic solutions (i.e., periodic solutions of higher period in the variable $\theta$) are fixed points of the iterated map~$\mathcal{P}_\varepsilon^m$ for some $m>1$, but not of $\mathcal{P}_\varepsilon$ itself. Furthermore, the quasi-periodic curves of the map~$\mathcal{P}_\varepsilon$ correspond to quasi-periodic solutions of \eqref{Continuous0}, each representing a dense orbit on an invariant torus.

\section{Proof of Theorem~\ref{MainTheorem}: non-ergodicity}\label{Sec4-Nonerg}

In this section we prove the first part of Theorem~\ref{MainTheorem}, which follows from an easy application of Moser's twist theorem to the map \eqref{Discrete0}. We conclude that every level set $\Scal_h$ of the Hamiltonian~\eqref{Ham} for $h$ small enough, exhibits a family of invariant tori whose union has positive measure. In particular, the associated Hamiltonian flow is non-ergodic for small enough energies. To prove this, we only use conditions \eqref{nonresonant} and \eqref{Isochronic}, which are both open and dense in $C^\infty(\TT^2)$, thus implying that generic potential fields are non-ergodic.\\
\begin{The}\label{Theorem-Twist}
	Let $V\in\mathcal{V}_0$ be a potential satisfying \eqref{nonresonant} and \eqref{Isochronic} (open and dense conditions), and $\varepsilon>0$ a small enough constant. Then, on every annulus $\mathcal{A}\times\TT\subset(0,d_I)\times\TT$, there is a positive measure set of smooth and closed invariant curves $\zeta$ of $\mathcal{P}_\varepsilon$. In addition, each $\zeta\subset\mathcal{A}\times\TT$ can be parameterized as
	\begin{equation*}
		I=I_0+U(\Psi),\quad\varphi=\varphi_0+V(\Psi),\quad\Psi\in\TT,
	\end{equation*}
	where $U,V$ are $\Ccal^1$, $2\pi$-periodic and 
	\begin{equation*}
		\|U\|_{\Ccal^1}+\|V\|_{\Ccal^1}<\varepsilon.
	\end{equation*} Moreover, the induced map $\mathcal{P}_\varepsilon:\zeta\rightarrow\zeta$ is given by 
	\begin{equation}\label{TorusParametrization}
		\phi\mapsto\phi+2\pi\sigma,
	\end{equation}for a certain irrational number $\sigma$ that satisfies the Diophantine condition
	\begin{equation}\label{diophantine}
		\left|\sigma-\frac{N}{M}\right|\geq\tau M^{-m},
	\end{equation}
	for some $\tau,m>0$, and all $M,N\in\mathbb{N}$.
\end{The}
\proof Take an arbitrary $\theta_0$ (the twist condition~\eqref{Isochronic} does not depend on this), and recall, cf. Remark~\ref{Twist}, that $|\Upsilon'(I)|$ only depends on the values $\{\omega_1,\omega_2, \beta_1,\beta_2,\beta_3\}$, which are fixed for each $V\in\mathcal{V}_0$.

On the other hand, the functions $f$ and $g$ defined in~\eqref{Discrete0} are integrals of finite sums of homogeneous polynomials, and hence they are uniformly bounded on any annulus $\mathcal{A}\times\TT\subset(0,d_I)\times\TT$. Therefore, for any $r\in\mathbb{N}$, there exists a positive constant $C$, depending on the $C^r$-norm of the Hamiltonian \eqref{BirkNF} on $\mathcal{A}\times\TT$, such that
\begin{equation}
	\|f\|_{C^r}+\|g\|_{C^r}<\varepsilon^3 |\Lambda_V|C.
\end{equation}
In this setting, Moser's twist theorem \cite[Theorem 4.8.2.]{GH} applies for every $\varepsilon\in(0,\bar\varepsilon)$, implying that the map $\mathcal{P}_\varepsilon$ has an invariant curve on $\mathcal{A}\times\TT$ associated with any irrational number $\sigma\in\text{Image}(\Upsilon)$ satisfying the Diophantine condition~\eqref{diophantine}. It is well known that the set of irrational numbers satisfying~\eqref{diophantine} with a fixed $m>2$ has full measure in every real interval. This completes the proof.
\cqfd

\section{Proof of Theorem~\ref{MainTheorem}: hyperbolic periodic orbits}\label{S:hyper}

Finally, in this section we complete the proof of Theorem~\ref{MainTheorem}. To this end, we analyze the dynamics near the resonant tori and prove that a generic potential $V$ exhibits an arbitrarily high number of hyperbolic periodic orbits on each small energy level. In particular, we prove the following result:
\begin{The}\label{TheoremHyperbolic}
	Let $N_0$ be a positive integer. There exists a generic set $\mathcal V'\subset C^\infty(\TT^2)$ and a small constant $\delta_0$ such that the Hamiltonian system $H_V$ exhibits at least $N_0$ elliptic and $N_0$ hyperbolic periodic orbits on each energy level $H^{-1}_V(h)$ for $h\in(0,\delta_0)$.
\end{The}

It is clear that Theorem~\ref{MainTheorem} is a straightforward consequence of this result and Theorem~\ref{Theorem-Twist}. Indeed, the set of potentials $\mathcal V'$ is residual, and the set of ``non-ergodic'' potentials given by Theorem~\ref{Theorem-Twist} is open and dense, so their intersection is a residual set as well. This is the set $\mathcal V\subset C^\infty(\mathbb T^2)$ in the statement of Theorem~\ref{MainTheorem}. 

To prove Theorem~\ref{TheoremHyperbolic}, we fix an integer $N_0\geq2$ and define the class $\mathcal{V}_{N_0}\subset\mathcal{V}_0$ of potentials satisfying \eqref{nonresonant}-\eqref{Isochronic}, and such that $\tfrac{\omega_1}{\omega_2}=\tfrac{n}{m}\in\mathbb{Q}$, where
\begin{equation}\label{coprimes}
		n\geq2,\ m\geq N_0\hbox{ are coprime and such that }[n+m]_{\ZZ_5}\leq2.
\end{equation} 
We prove in Appendix~\ref{Appendix-Denseclass} that the set $\mathcal{V}_{N_0}$ is dense in $ C^\infty(\TT^2)$. Given $V\in\mathcal{V}_{N_0}$, consider the symplectic transformation introduced in Corollary~\ref{Cor-SymplecticTransformation} and set $\theta_0=0$, without any loss of generality. Then, this implies
\begin{equation*}
	\Upsilon(I)=nm^{-1}+\varepsilon^2\Lambda(I),
\end{equation*}
and, for any initial datum $(I(0),\varphi(0))=(I_0,\varphi_0)\in(0,d_I)\times\TT$, the $m$-iterated Poincaré map \eqref{DiscreteM} takes the form:
\begin{equation}\label{DiscreteM-Resonance}
	\left\{\begin{array}{l}
	I^m=I_0+f^m(I,\varphi,\varepsilon),\\
	 \varphi^m=\varphi_0+\varepsilon^2\int_{0}^{2\pi m}\Lambda(I(\theta))d\theta+g^m(I,\varphi,\varepsilon).
	\end{array}
	\right.
\end{equation}
In contrast with \eqref{DiscreteM}, we omit $\theta_0$ from the notation here, since it is fixed. Regarding the perturbative terms in \eqref{DiscreteM-Resonance}, we write them as follows:
\begin{align*}
	&f^m(I,\varphi,\varepsilon)=\sum\limits_{l=5}^{N}\varepsilon^{l-2}\int_{0}^{2\pi m}\dfrac{\partial_\varphi H_l(I(\theta),J(I(\theta)),\varphi(\theta),\theta)}{\partial_J H_V(I(\theta),J(I(\theta)),\varphi(\theta),\theta)}d\theta+\Theta(\varepsilon^{N+1}),\\
	&g^m(I,\varphi,\varepsilon)=\sum\limits_{l=5}^{N}\varepsilon^{l-2}\int_{0}^{2\pi m}\dfrac{\partial_I H_l(I(\theta),J(I(\theta)),\varphi(\theta),\theta)}{\partial_J H_V(I(\theta),J(I(\theta)),\varphi(\theta),\theta)}d\theta+\Theta(\varepsilon^{N+1}),
\end{align*}
with arbitrary $N\gg N_0$. We recall that any fixed point of \eqref{DiscreteM-Resonance} corresponds to a $2\pi m$-periodic orbit of \eqref{Continuous0}. Next, we aim to study the existence of fixed points of the $m$-iterated Poincar\'e map.

First, we focus on studying the expressions:
 \begin{equation}\label{Action-Force}
 	\sum\limits_{l=3}^{N}\varepsilon^{l}\dfrac{\partial_\nu H_{l+2}}{\partial_J H_V}(I,J,\vec{\varphi})+\Theta(\varepsilon^{N+1}),\ \hbox{ for }\nu\in\{\varphi,I\},
 \end{equation}
{for arbitrary $N\geq3$}, by describing explicitly their dependence on the perturbation parameter $\varepsilon$, which is strongly nonlinear. Concerning this, let us recall from \eqref{levelSet} that the reversed action $J$ is a function such that
\begin{equation}\label{J-function}
	J=J(I,\vec{\varphi},\varepsilon;h)=J_0(I)+\sum_{l=2}^{N}\varepsilon^l J_l(I,\vec{\varphi};h)+\Theta(\varepsilon^{N+1}),
\end{equation}
for any $N\in\mathbb{N}$, and for all sufficiently small $h>0$. In Section~\ref{Section-ReversedAction} we study the structure of $J$ in detail. For this purpose, for any $k\in\{5,6,7\}$ and {for any sufficiently large} $\mm\in\mathbb{N}$, we define the sets
\begin{equation}\label{setL0}
	\mathcal{L}_0^k[1,\mm]:=\left\{H_k,\ \partial_J^{s} H_k:\ 1\leq s\leq\mm\right\}.
\end{equation}
{Similarly, for $\nu\in\{\varphi,I\}$ we define}
\begin{equation*}
\mathcal{L}_{0,\nu}^k[1,\mm]:=\mathcal{L}_0^k[1,\mm]\cup\left\{\partial_\nu H_k,\ \partial^{J}_\nu H_k,\ \partial^{J,s}_\nu H_k:\ 2\leq s\leq \mm\right\},
\end{equation*}where we write
\begin{equation*}
	\partial^{J}_\nu H_k:=\partial_J\partial_\nu H_k,\quad\hbox{and}\quad\partial^{J,s}_\nu H_k:=\partial^s_J\partial_\nu H_k,
\end{equation*} to simplify the notation. Moreover, all functions are evaluated at $(I,J_0(I),\vec{\varphi})$.

Similarly, for any $l\geq 2$, we define $\mathcal{L}_0^k[l,\mm]$ as the set of finite linear combinations of $l$-products of elements from $\mathcal{L}_0^k[1,\mm]$. Since $\mathcal{L}_{0}^k[1,\mm]\subset\mathcal{L}_{0,\nu}^k[1,\mm]$, it follows that $\mathcal{L}_0^k[l,\mm]$ is contained in the set of finite linear combinations of $l$-products of elements from $\mathcal{L}_{0,\nu}^k[1,\mm]$, for both $\nu\in\{\varphi,I\}$. Accordingly, we define $\mathcal{L}_{0,\nu}^k[l,\mm]$ as the relative complement of $\mathcal{L}_0^k[l,\mm]$ in each of these larger sets.\\ Additionally, for $k\in\{6,7\}$, we say that $\ff\in\mathcal{L}_1^k[l,\mm]$ if it is a finite linear combination of the form
\begin{equation*}
	\ff=\sum_{j}\ff_j\g_j,\quad\hbox{with } (\ff_j,\g_j)\in\mathcal{L}_0^5[l-1,\mm]\times\mathcal{L}_0^k[1,\mm],\hbox{ for all }j.
\end{equation*}
Furthermore, let $\mathcal{L}^7_2[l,\mm]$ be the set of functions $\ff$ with the following structure:
\begin{equation*}
\ff=\sum_{j}\left[\ff_{j}\g_j+\mathbf{h}_j\right],\quad\hbox{with } (\ff_j,\g_j,\mathbf{h}_j)\in\mathcal{L}_0^5[l-2,\mm]\times\mathcal{L}_0^6[2,\mm]\times\mathcal{L}_1^7[l,\mm],\hbox{ for all }j.
\end{equation*}
Finally, for both $\nu\in\{\varphi,I\}$, the spaces $\mathcal{L}^7_{2,\nu}$ and $\mathcal{L}^k_{1,\nu}$ (with $k\in\{6,7\}$) are defined analogously to $\mathcal{L}^k_{0,\nu}$.
\begin{Rem}\label{Rem-ResonanceCondition}
 Since $J_0$ becomes constant when $I\in(0,d_I)$ is fixed, any function in the previously defined spaces is, in that case, a homogeneous polynomial of a certain degree with respect to the angular variables
\begin{equation}\label{levelSet-AngularVariables}
	\xi_1=e^{-i\varphi},\quad \xi_2=e^{-i\theta},\quad \eta_k=i\bar\xi_k,\quad k=1,2.
\end{equation} 
\end{Rem}
On the other hand, for any $\gamma\in\Gamma_k$, with $\Gamma_k$ defined in \eqref{set-Gamma_l} and $k\in\{5,6,7\}$, let $\mathbf{n}_\gamma$ be the following number:
\begin{equation*}
\mathbf{n}_\gamma=\left\{\begin{array}{ll}
	\tfrac{\gamma_2+\gamma_4}{2},&\hbox{if $\gamma_2+\gamma_4$ is even,}\\
	k,&\hbox{ otherwise,}
\end{array}\right.
\end{equation*}
and consider the quantities:
\begin{equation*}
	C_I(\gamma,0):=\sqrt{I^{\gamma_1+\gamma_3}J_0(I)^{\gamma_2+\gamma_4}},
\end{equation*} and
\begin{equation}\label{Constant-I}
C_I(\gamma,s):=\prod_{j=1}^{s}\left(\tfrac{\gamma_2+\gamma_4}{2}-j+1\right)I^{\tfrac{\gamma_1+\gamma_3}{2}}J_0(I)^{\tfrac{\gamma_2+\gamma_4}{2}-s},\ \hbox{for }1\leq s\leq\mathbf{n}_\gamma,
\end{equation}
with $s\in\mathbb N$. Then, we define the set
\begin{equation*}
\Pi_0^{k}:=\left\{C_\gamma C_I(\gamma,s) :\ \gamma\in\Gamma_k,\ 0\leq s\leq \mathbf{n}_\gamma \right\},
\end{equation*} where the constant $C_\gamma$ is defined in \eqref{Constant-FinalForm}. Moreover, let $\Pi_0^{k,l}$ denote the set of all finite linear combinations of $l$-products of elements from  $\Pi_0^{k}$. Similarly,
\begin{align}\label{ConstantsSet}
\Pi_{0,\varphi}^{k}:=\{C_\gamma C_I(\gamma,s),\ i(\gamma_3-\gamma_1)C_\gamma C_I(\gamma,s):\ \gamma\in\Gamma_k,\ 0\leq s\leq \mathbf{n}_\gamma \},
\end{align}
and 
\begin{align}\label{ConstantsSet'}
	\Pi_{0,I}^{k}:=\{C_\gamma C_I(\gamma,s),\ C_\gamma \partial_IC_I(\gamma,s):\ \gamma\in\Gamma_k,\ 0\leq s\leq \mathbf{n}_\gamma \}.
\end{align}
In addition, for both $\nu\in\{\varphi,I\}$, we define $\Pi_{0,\nu}^{k,l}$ as the relative complement of $\Pi_{0}^{k,l}$ in the set of all finite linear combinations of $l$-products of elements from $\Pi_{0,\nu}^{k}$. As a result, any function $\ff\in\mathcal{L}^k_{0,\nu}[l,\mm]$ (or in $\mathcal{L}^k_{0}[l,\mm]$) can be written as follows:
\begin{equation}\label{J-Constants}
\ff=\sum\limits_{\gamma\in\Gamma_{kl}} i^{\gamma_3+\gamma_4}\mathbf{C}[\gamma,I]e^{-i\varphi(\gamma_1-\gamma_3)}e^{-i\theta(\gamma_2-\gamma_4)},
\end{equation}
with $\mathbf{C}[\gamma,I]\in\Pi_{0,\nu}^{k,l}$ (or in $\Pi_{0}^{k,l}$), for all $\gamma\in\Gamma_{kl}$.  
\begin{Rem}\label{Rem-OpenDenseConstant}
Observe that the constants $C_I$ and $\partial_I C_I$ are all finite and nonzero when $I\neq0$. Consequently, for any fixed pair $(\gamma,I)\in\Gamma_{kl}\times(0,d_I)$, the condition
\begin{equation*}
	\mathbf{C}[\gamma,I]\neq0
\end{equation*}
is open and dense in $ C^\infty(\TT^2)$.
\end{Rem} 
On the other hand, it is natural to define the sets $\Pi_{1}^{k,l}$ (with $k\in\{6,7\}$) and $\Pi_2^{7,l}$ as follows:
\begin{equation*}
C\in\Pi_{1}^{k,l}\Leftrightarrow C=\sum_{j} \mathbf{a}_j\mathbf{b}_j,\hbox{ with }(\mathbf{a}_j,\mathbf{b}_j)\in\Pi_0^{5,l-1}\times\Pi_0^{k},\hbox{ for all }j,
\end{equation*}
and
\begin{equation*}
C\in\Pi_2^{7,l}\Leftrightarrow C=\sum_{j}[ \mathbf{a}_j\mathbf{b}_j+\mathbf{c}_j],\hbox{ with }(\mathbf{a}_j,\mathbf{b}_j,\mathbf{c}_j)\in\Pi_0^{5,l-2}\times\Pi_0^{6,2}\times\Pi_1^{7,l},\hbox{ for all }j,
\end{equation*}
while $\Pi_{1,\nu}^{k,l}$ and $\Pi_{2,\nu}^{7,l}$ are constructed analogously to $\Pi_{0,\nu}^{k,l}$, for both $\nu\in\{\varphi,I\}$.

To conclude, observe that any element of $\mathcal{L}_{1}^{k,l}$, $\mathcal{L}_{1,\nu}^{k,l}$, $\mathcal{L}_{2}^{7,l}$, or $\mathcal{L}_{2,\nu}^{7,l}$, is also a polynomial of the form given in \eqref{J-Constants}, but with degree $5l+[k]_{\ZZ_5}$, and with the constants $\mathbf{C}[\gamma, I]$ all belonging to the corresponding space $\Pi_{1}^{k,l}$, $\Pi_{1,\nu}^{k,l}$, $\Pi_{2}^{7,l}$, or $\Pi_{2,\nu}^{7,l}$. Therefore, it is clear that Remark~\ref{Rem-OpenDenseConstant} is also valid in these cases.\\

Now we state a result that describes explicitly the dependence on $\varepsilon$ of the functions \eqref{Action-Force}, and whose proof is given in Section~\ref{Section-Proof1}. 
\begin{Pro}\label{Pro-Hyp0}
	Let $V\in\mathcal{V}_0$ be a potential satisfying \eqref{nonresonant}. Then, on every level set \eqref{levelSet} with sufficiently small energy $h>0$, it holds that
	\begin{align*}
		\sum\limits_{l=3}^{N}\varepsilon^{l}\dfrac{\partial_\varphi H_{l+2}(I,J,\vec{\varphi})}{\partial_J H_V(I,J,\vec{\varphi})}=\sum\limits_{l=3}^{N}\varepsilon^l F_l(I,J_0(I),\vec{\varphi})+\Theta(\varepsilon^{N+1}),\\
		\sum\limits_{l=3}^{N}\varepsilon^{l}\dfrac{\partial_I H_{l+2}(I,J,\vec{\varphi})}{\partial_J H_V(I,J,\vec{\varphi})}=\sum\limits_{l=3}^{N}\varepsilon^l G_l(I,J_0(I),\vec{\varphi})+\Theta(\varepsilon^{N+1}),
	\end{align*}
	{for any $N\geq3$}, and where every $F_l$, $G_l$ is a polynomial of degree $l+2\left\lfloor\tfrac{l}{3}\right\rfloor$ with respect to the angular variables \eqref{levelSet-AngularVariables}. More concretely,
		\begin{equation}
			F_l=\left(\dfrac{-1}{\omega_2}\right)^{\left\lfloor\tfrac{l}{3}\right\rfloor}\ff_l+\hbox{l.o.t;}\quad 	G_l=\left(\dfrac{-1}{\omega_2}\right)^{\left\lfloor\tfrac{l}{3}\right\rfloor}\g_l+\hbox{l.o.t.}
		\end{equation} 
		for some explicit $\ff_l\in\mathcal{L}_{[l]_{\ZZ_3},\varphi}^{5+[l]_{\ZZ_3}}\left[\left\lfloor\tfrac{l}{3}\right\rfloor,\left\lfloor\tfrac{l}{3}\right\rfloor-1\right],$ and $\g_l\in\mathcal{L}_{[l]_{\ZZ_3},I}^{5+[l]_{\ZZ_3}}\left[\left\lfloor\tfrac{l}{3}\right\rfloor,\left\lfloor\tfrac{l}{3}\right\rfloor-1\right]$.
	
\end{Pro}

On the other hand, given a pair $\{n,m\}$ satisfying \eqref{coprimes}, there is no vector $\gamma\in\mathbb{N}^4$ with $|\gamma|<n+m$ and such that
\begin{equation}\label{ResonanceEquality}
	n(\gamma_1-\gamma_3)=m(\gamma_4-\gamma_2)\neq0.
\end{equation}
When $|\gamma|=n+m,$ then \eqref{ResonanceEquality} is satisfied only for $$\gamma\in\{\gamma^m_1,\gamma^m_2\}:=\{(m,0,0,n),\ (0,n,m,0)\}.$$
Moreover, we observe that there exists a unique $l\geq3$, which we denote by $\LL$, such that $$m+n=\LL+2\left\lfloor\tfrac{\LL}{3}\right\rfloor\,.$$  
Then, by combining these observations, Remark~\ref{Rem-ResonanceCondition} and Proposition~\ref{Pro-Hyp0}, we obtain the following result:
\begin{Cor}\label{Cor-Hyp}
Given $V\in\mathcal{V}_{N_0}$, consider a level set \eqref{levelSet} with sufficiently small energy $h>0$, and fix an arbitrary $(I_0,\varphi_0)\in(0,d_I)\times\TT$. Then, 
\begin{equation*}
	\int_{0}^{2\pi m}F_l\left(\hat{X}(\theta)\right)d\theta=\int_{0}^{2\pi m}G_l\left(\hat{X}(\theta)\right)d\theta=0,\quad\hbox{ for all }l<\LL,
\end{equation*}
where 
\begin{equation}\label{XTheta}
	\hat X(\theta):=(I_0,J_0(I_0),\varphi_0+nm^{-1}\theta,\theta).
\end{equation}
Moreover,
\begin{align}\label{Melnikov1}
	\int_{0}^{2\pi m}F_{\LL}\left(\hat{X}(\theta)\right)d\theta=&\cos(m\varphi_0)\left[i^m\mathbf{C}_{\varphi}[\gamma^m_2,I_0]+i^n\mathbf{C}_\varphi[\gamma^m_1,I_0]\right]\\&-i\sin(m\varphi_0)\left[i^m\mathbf{C}_{\varphi}[\gamma^m_2,I_0]-i^n\mathbf{C}_{\varphi}[\gamma^m_1,I_0]\right],
\end{align}
and
\begin{align}\label{Melnikov1'}
	\int_{0}^{2\pi m}G_{\LL}\left(\hat{X}(\theta)\right)d\theta=&\cos(m\varphi_0)\left[i^m\mathbf{C}_I[\gamma^m_2,I_0]+i^n\mathbf{C}_I[\gamma^m_1,I_0]\right]\\&-i\sin(m\varphi_0)\left[i^m\mathbf{C}_I[\gamma^m_2,I_0]-i^n\mathbf{C}_I[\gamma^m_1,I_0]\right].
\end{align}
\end{Cor}
\begin{Rem}\label{RemarkMelnikov}
	For any $V\in\mathcal{V}_0$, both functions $F_l$ and $G_l$ are real-valued. Therefore, the integrals
	\begin{equation*}
		\int_{0}^{2\pi m}F_l\left(\hat{X}(\theta)\right)d\theta,\quad\hbox{and}\quad \int_{0}^{2\pi m}G_l\left(\hat{X}(\theta)\right)d\theta,
	\end{equation*}also take real values, for all $l\geq 3$, and all $(I_0,\varphi_0)\in(0,d_I)\times\TT$. Consequently, the constants in \eqref{Melnikov1}-\eqref{Melnikov1'} are characterized as follows:
	\begin{equation*}
		i^{n}\mathbf{C}_\nu[\gamma^m_1,I_0]=-i^m\overline{\mathbf{C}_\nu[\gamma^m_2,I_0]},\quad\hbox{for }\nu\in\{I,\varphi\},
	\end{equation*} 
	so that \eqref{Melnikov1}-\eqref{Melnikov1'} can also be written as
	\begin{equation*}
	2\left(\cos(m\varphi_0)\textnormal{Re}\left(i^n\mathbf{C}_\nu[\gamma^m_1,I_0]\right)-\sin(m\varphi_0)\textnormal{Im}\left(i^n\mathbf{C}_\nu[\gamma^m_1,I_0]\right)\right),\ \hbox{for }\nu\in\{I,\varphi\},
    \end{equation*}
	where, as usual, $\textnormal{Re}(x)$ and $\textnormal{Im}(x)$ denote the real and imaginary parts of $x\in\mathbb{C}$, respectively.
\end{Rem}
Finally, the next result provides an explicit expression for the $m$-iterated Poincaré map introduced in Equation~\eqref{DiscreteM-Resonance}. The integer $\LL$ introduced before plays a key role in these formulas. The proof of this proposition is presented in Section~\ref{Section-Proof2}.
\begin{Pro}\label{Pro-Hyp1}
	Given $V\in\mathcal{V}_{N_0}$, consider a level set \eqref{levelSet} of sufficiently small energy $h>0$. Then, for any initial condition $(I_0,\varphi_0)\in(0,d_I)\times\TT$ of the Poincaré map \eqref{DiscreteM-Resonance}, the following holds:
\begin{align*}
\varepsilon^2\int_0^{2\pi m}\Lambda(I(\theta))d\theta=\varepsilon^2 2\pi m\Lambda(I_0)+\Theta(\varepsilon^{\LL+2}),
\end{align*}
and
\begin{equation}\label{Melnikov2}	
\left.\begin{array}{l}
f^m(I,\varphi,\varepsilon)=2\varepsilon^{\LL}\Big(\cos(m\varphi_0)\textnormal{Re}(\textnormal{\textbf{C}}^\varphi)-\sin(m\varphi_0)\textnormal{Im}({\textbf{C}}^\varphi)\Big)+\Theta(\varepsilon^{\LL+1}),
\vspace{1mm}\\
g^m(I,\varphi,\varepsilon)=2\varepsilon^{\LL}\left(\cos(m\varphi_0)\textnormal{Re}(\textnormal{\textbf{C}}^I)-\sin(m\varphi_0)\textnormal{Im}(\textnormal{\textbf{C}}^I)\right)+\Theta(\varepsilon^{\LL+1}),
\end{array}\right.
\end{equation}
for some explicit constants $\left\{\textnormal{\textbf{C}}^\nu: \ \nu\in\{I,\varphi\}\right\}$.
\end{Pro}
\begin{Rem}
From the proof of Proposition~\ref{Pro-Hyp1} it follows that 
\begin{equation}
{\mathbf{C}}^\nu=i^n\mathbf{C}_\nu[\gamma^m_1,I_0]+i^n\mathbf{C}_{\nu,r}[\gamma^m_1,I_0],\ \hbox{ for }\nu\in\{I,\varphi\},
\end{equation}
where the constant $\mathbf{C}_{\nu}[\gamma^m_1,I_0]\in\Pi^{5+[m+n]_{\ZZ_3}}_{[m+n]_{\ZZ_3},\nu}$, cf.  Remark~\ref{RemarkMelnikov}. { Moreover, for both $\nu\in\{I,\varphi\},$ the remainder $\mathbf{C}_{\nu,r}[\gamma^m_1,I_0]$ is an algebraic function of elements from \eqref{ConstantsSet} and \eqref{ConstantsSet'}. In fact, $\mathbf{C}_{\nu,r}[\gamma^m_1,I_0]$ necessarily involves elements from both sets and cannot depend solely on elements from one of them. On the contrary, $\mathbf{C}_{\nu}[\gamma^m_1,I_0]$ depends solely on elements from $\bigcup\limits_{k}\Pi^k_{0,\nu}.$ }Therefore, for any fixed $I_0\in(0,d_I)$, the condition
\begin{equation}\label{MelnikovOpenDenseCondition}
{\mathbf{C}}^\nu\neq0,
\end{equation} is open and dense in $ C^\infty(\TT^2)$, for both $\nu\in\{I,\varphi\}$.
\end{Rem}

Using all these results we are ready to complete the proof of Theorem~\ref{TheoremHyperbolic}. Indeed, for any initial datum $(I_0,\varphi_0)\in(0,d_I)\times\TT$, combining Proposition~\ref{Pro-Hyp1}, Corollary~\ref{Cor-Hyp}, and Remark~\ref{RemarkMelnikov}, the $m$-iterated Poincaré map~\eqref{DiscreteM-Resonance} reads as:
\begin{equation}\label{DiscreteM-Resonance2}
	\left\{\begin{array}{cll} 
	\dfrac{I^m-I_0}{2\varepsilon^\LL}&=&f(I_0,\varphi_0)+\Theta(\varepsilon),\vspace{1mm}\\
	 \dfrac{\varphi^m-\varphi_0}{m2\pi\varepsilon^2}&=&\Lambda(I_0)+\Theta(\varepsilon^{\LL-2}),
	\end{array}
	\right.
\end{equation}
where $f:(0,d_I)\times\TT\to\RR$ is the trigonometric function:
\begin{equation}\label{trigonometric}
	f(I,\varphi):=\cos(m\varphi)\textnormal{Re}\left(\Cc^\varphi\right)-\sin(m\varphi)\textnormal{Im}\left(\Cc^\varphi\right).
\end{equation}
Therefore, when $\Cc^\varphi\neq0$,
the function \eqref{trigonometric} has $2m$ roots $\{\varphi_0^k:\ 1\leq k\leq 2m\}$ in $\TT$, all of which are nondegenerate. Accordingly, if we assume the generic condition \eqref{MelnikovOpenDenseCondition} and we fix $I_0=I_0^*$ as introduced in Equation~\eqref{ResonanceCondition}, it follows that the pairs $\{(I_0^*,\varphi_0^k):\ 1\leq k\leq 2m\}$ are the zeros of the non-perturbative part in Equation~\eqref{DiscreteM-Resonance2}. Furthermore, the Jacobian at each of them is given by
\begin{equation}\label{JacobianPoincare}
\mathcal{J}(I_0^*,\varphi_0^k)=m\Lambda_V\left(\sin(m\varphi_0^k)\textnormal{Re}\left(\Cc^\varphi\right)+\cos(m\varphi_0^k)\textnormal{Im}\left(\Cc^\varphi\right)\right)\neq0,
\end{equation}
for all $1\leq k\leq 2m$. It then follows from the Implicit Function Theorem that there exists a value $\varepsilon^*>0$ and $2m$ smooth functions
\begin{equation*}
	(I_\varepsilon^k,\varphi_\varepsilon^k):[0,\varepsilon^*)\rightarrow(0,d_I)\times\TT,\hbox{ with }1\leq k\leq2m,
\end{equation*}
which are the zeros of the map \eqref{DiscreteM-Resonance2}, and such that 
\begin{equation}\label{JacobianPoincare2}
	\mathcal{J}(I_\varepsilon^k,\varphi_\varepsilon^k)=\mathcal{J}(I_0^*,\varphi_0^k)+\Theta(\varepsilon)\neq0,
\end{equation}
for all $\varepsilon\in[0,\varepsilon^*)$ and all $1\leq k\leq 2m$. Equivalently, for each $\varepsilon\in(0,\varepsilon^*)$, the set $\{(I_\varepsilon^k,\varphi_\varepsilon^k): 1\leq k\leq 2m\}$ is the set of fixed points of the Poincaré map \eqref{DiscreteM-Resonance}, and
\begin{equation*}
(I_\varepsilon^k,\varphi^k_\varepsilon)\rightarrow(I_0^*,\varphi^k_0)\ \hbox{ when }\varepsilon\rightarrow0,\hbox{ for all }1\leq k\leq 2m.
\end{equation*}
Next we compute the eigenvalues $(\lambda_+^k,\lambda_-^k)$ of the map \eqref{DiscreteM-Resonance} at the fixed points $\{(I_\varepsilon^k,\varphi_\varepsilon^k)\}$ in order to identify their stability type. By Proposition~\ref{Pro-Hyp1}, the corresponding characteristic polynomials are of the form:
\begin{align*}
0=&(1-\lambda)^2+(1-\lambda)\left[\partial_I f^m(I_\varepsilon^k,\varphi_\varepsilon^k,\varepsilon)+\partial_\varphi g^m(I_\varepsilon^k,\varphi_\varepsilon^k,\varepsilon)\right]\\
&+\left[\varepsilon^{2}m2\pi\Lambda_V+\partial_Ig^m(I_\varepsilon^k,\varphi_\varepsilon^k,\varepsilon)\right]\partial_\varphi f^m(I_\varepsilon^k,\varphi_\varepsilon^k,\varepsilon)+\Theta(\varepsilon^{\LL+2})\\
=&(1-\lambda)^2+(1-\lambda)\Theta(\varepsilon^{\LL})+\varepsilon^{\LL+2}4m\pi\mathcal{J}(I_\varepsilon^k,\varphi_\varepsilon^k)+\Theta(\varepsilon^{\min\{2\LL,\LL+4\}}),
\end{align*}
which implies that
\begin{equation}\label{eigenvalues-hyperbolic}
	1-\lambda_{\pm}=\pm\varepsilon^{(\LL+2)/2}\sqrt{4m\pi\mathcal{J}(I_\varepsilon^k,\varphi_\varepsilon^k)}+\Theta(\varepsilon^{\min\{\LL,\LL/2+2\}}).
\end{equation}
Observe that the Jacobian \eqref{JacobianPoincare} is positive at $m$ of the pairs $\{I_0^*,\varphi_0^k\}$ and negative at the remaining ones. Since the same holds for \eqref{JacobianPoincare2}, it follows that half of the fixed points are elliptic, while the other half are hyperbolic. Moreover, the set of potentials $V\in\mathcal{V}_{N_0}$ that satisfy the (open and dense) condition~\eqref{MelnikovOpenDenseCondition} is also dense in $C^\infty(\TT^2)$, and we denote such a set of potentials by $\widetilde{\mathcal V}$.

To conclude, we {recall that, by construction, $m\geq N_0$ for all $V\in\widetilde{\mathcal V}$, and that the computation of the fixed points works for every $\varepsilon\in(0,\varepsilon^*)$. Equivalently, in view of Remark~\ref{Remark-EnergyLevelSets}, the result holds for all energy values 
\begin{equation}
    0<h<(\varepsilon^*)^2h_0=:\delta_0, 
\end{equation}
for some positive $h_0$, which is a small ($\varepsilon$-independent) parameter in order to satisfy both the level set equation \eqref{levelSet} and condition \eqref{nonresonant}. In particular, $\delta_0\ll1$. Therefore, since fixed points of the Poincar\'e map correspond to periodic orbits (with the same stability type) of the Hamiltonian system, the Hamiltonian system $H_V$ exhibits at least $N_0$ hyperbolic periodic orbits on each energy level $H^{-1}_V(h)$ with $h\in(0,\delta_0)$, and $V\in\widetilde{\mathcal V}$. } Covering the open interval $(0,\delta_0)$ with countably many compact intervals $\{K_j\}_{j=-\infty}^\infty$, the hyperbolic permanence theorem implies that for each $V\in \widetilde{\mathcal V}$ there is an open neighborhood of potentials which also exhibit $N_0$ hyperbolic periodic orbits on each energy level $H^{-1}_V(h)$ with $h\in K_j$. This clearly yields an open and dense set of potentials for each $j\in\ZZ$. Theorem~\ref{TheoremHyperbolic} then follows by taking the countable intersection $\mathcal V'$ of these open and dense sets, which is residual because $C^\infty(\TT^2)$ is Baire.

\section{Structure of the reversed action function}\label{Section-ReversedAction}
In this section, we analyze the reversed action function $J$, defined by \eqref{J-function} on energy levels close to the minimum. We recall that, on these level sets of small energy $h\ll1$, $J:(0,d_I)\times\TT^2\times[0,\varepsilon^*)\to(0,d_J)$ is a well-defined and smooth function with the following structure:
\begin{equation*}
	J(I,\vec{\varphi},\varepsilon;h)=J_0(I)+\sum_{l=2}^{N}\varepsilon^l J_l(I,\vec{\varphi};h)+\Theta(\varepsilon^{N+1}),
\end{equation*}
where the choice of $N\in\mathbb{N}$ is arbitrary. The following proposition is the main result of this section:
\begin{Pro}\label{J-Proposition} Let $V\in\mathcal{V}_0$ be satisfying \eqref{nonresonant}. Then, on every level set \eqref{levelSet} with sufficiently small energy $h>0$, every function $J_l$ is a polynomial of degree $l+2\lfloor\tfrac{l}{3}\rfloor$ in the angular variables defined in \eqref{levelSet-AngularVariables}. More concretely, 
	\begin{equation}
		J_l=\left(\dfrac{-1}{\omega_2}\right)^{\left\lfloor\tfrac{l}{3}\right\rfloor}\jj_l\ +\hbox{ l.o.t.}
	\end{equation}
	for some explicit $\jj_l\in\mathcal{L}_{[l]_{\ZZ_3}}^{5+[l]_{\ZZ_3}}\left[\left\lfloor\tfrac{l}{3}\right\rfloor,\left\lfloor\tfrac{l}{3}\right\rfloor-1\right].$ 
\end{Pro} 

The proof of this proposition is divided in two subsections. First, we establish several key lemmas, which are used then to prove the result.

\subsection{Auxiliary lemmas}
Most of the computations in our analysis rely on the application of the Faà di Bruno's formula of the generalized chain rule~\cite[1.13]{KP}. Concretely, we shall use it to compute higher order derivatives in Equations~\eqref{levelSet}-\eqref{Continuous1}, for which it is convenient to define the sets 
\begin{equation}\label{SetGamma}
	\widetilde{\Gamma}_k=\{\gamma\in\mathbb{N}^{k-3}:\sum\limits_{j} (j+1)\gamma_j=k\},\ \ \hbox{with }k\geq4.
\end{equation}
Note that the definition of $\widetilde{\Gamma}_k$ is linked to the fact that $\left.\partial_\varepsilon J\right|_{\varepsilon=0}=J_1=0$.

On the one hand, by substituting $J$ at $H_1(I,J)$ we obtain
\begin{align}
	H_1(I,J)= H_1(I,J_0)+\dfrac{\beta_2}{2}\left[\sum\limits_{k=2}^N\varepsilon^{k} J_k\right]^2+\partial_{J}H_1(I,J_0)\sum\limits_{k=2}^N\varepsilon^{k} J_k+\Theta(\varepsilon^{N+1}).
\end{align}
Additionally, for any function in $\{H_l(I,J,\vec{\varphi}): l\geq5 \}$, its $N^{\text{th}}$-order Taylor expansion at $\varepsilon=0$ is
\begin{align}
	H_l(I,J,\vec{\varphi})=&\ H_l+\sum\limits_{k=2}^{N}\dfrac{\varepsilon^{k}}{k!}\partial_{J}H_l\partial_\varepsilon^{k}J+\sum\limits_{k=4}^{N}\dfrac{\varepsilon^{k}}{k!}\sum\limits_{\gamma\in\widetilde{\Gamma}_k}\dfrac{k!}{\gamma!}\partial^{|\gamma|}_JH_{l}\prod_{j=1}^{k-3}\left(\dfrac{\partial_\varepsilon^{j+1} J}{(j+1)!}\right)^{\gamma_{j}}+\Theta(\varepsilon^{N+1})\\=&\ H_l+\sum\limits_{k=2}^{N}\varepsilon^{k}\partial_{J}H_lJ_{k}+\sum\limits_{k=4}^{N}\varepsilon^{k}\sum\limits_{\gamma\in\widetilde{\Gamma}_k}\dfrac{\partial^{|\gamma|}_JH_{l}}{\gamma!}\prod_{j=1}^{k-3}J_{j+1}^{\gamma_{j}}+\Theta(\varepsilon^{N+1}),
\end{align}
where all terms are evaluated at $(I,J_0,\vec{\varphi})$. Therefore, these expressions reduce the level set equation \eqref{levelSet2} to
\begin{align}
	0=&\ \omega_2\sum\limits_{l=2}^{N}\varepsilon^l J_l+\varepsilon^2 H_1+\varepsilon^2\dfrac{\beta_2}{2}\left[\sum\limits_{l=2}^{N}\varepsilon^{l}J_l\right]^2+\varepsilon^2\partial_{J}H_1\sum\limits_{l=2}^{N}\varepsilon^{l} J_l\\
	&+\sum\limits_{l=5}^{N}\varepsilon^{l-2}\left[H_l+\sum\limits_{k=2}^{N}\varepsilon^{k}\partial_{J}H_lJ_{k}\right]+\sum\limits_{l=5}^{N}\sum\limits_{k=4}^{N}\varepsilon^{l+k-2}\sum\limits_{\gamma\in\widetilde{\Gamma}_k}\dfrac{\partial^{|\gamma|}_JH_{l}}{\gamma!}\prod_{j=1}^{k-3}J_{j+1}^{\gamma_{j}}+\Theta(\varepsilon^{N+1}),
\end{align}
and, after identifying the terms of order $\Theta(\varepsilon^{N+1})$, we can rearrange it as follows:
\begin{align}
	0=&\ \omega_2\sum\limits_{l=2}^{N}\varepsilon^l J_l+\varepsilon^2 H_1+\beta_2\left[\sum\limits_{l=3}^{\lfloor\tfrac{N}{2}\rfloor}\dfrac{\varepsilon^{2l}}{2}J^2_{l-1}+\sum_{l\geq3}^{\lfloor\tfrac{N}{2}\rfloor-1}\sum_{k\geq l+1}^{N-l}\varepsilon^{l+k}J_{l-1}J_{k-1}\right]\\&+\partial_{J}H_1\sum\limits_{l=4}^{N}\varepsilon^{l} J_{l-2}+\sum\limits_{l=3}^{N}\varepsilon^{l}H_{l+2}+\sum\limits_{l=5}^{N}\varepsilon^{l}\partial_{J}H_lJ_{2}+\sum\limits_{l=5}^{N-1}\sum\limits_{k=1}^{N-l}\varepsilon^{l+k}\partial_{J}H_lJ_{k+2}\\
	&+\sum\limits_{l=3}^{N-4}\sum\limits_{k=4}^{N-l}\varepsilon^{l+k}\sum\limits_{\gamma\in\widetilde{\Gamma}_k}\dfrac{\partial^{|\gamma|}_JH_{l+2}}{\gamma!}\prod_{j=1}^{k-3}J_{j+1}^{\gamma_{j}}+\Theta(\varepsilon^{N+1}).
\end{align}
 We can then compute each $J_l$ from the vanishing of each term of this power series in $\varepsilon$. Concretely,
\begin{equation}\label{J-function1}
	\left.\begin{array}{l}
		J_2=\dfrac{-H_1}{\omega_2},\qquad J_3=\dfrac{-H_5}{\omega_2},\qquad J_4=\dfrac{-H_6-\partial_J H_1J_2}{\omega_2},\vspace{1.2mm}\\
		J_5=\dfrac{-H_7-\partial_J H_1J_3-\partial_J H_5 J_2}{\omega_2},\vspace{0.5mm}\\
		J_6=\dfrac{-H_8-\partial_J H_1J_4-\partial_J H_6 J_2-\partial_JH_5J_3-\tfrac{\beta_2}{2}J^2_2}{\omega_2},
	\end{array}\right.
\end{equation}
and, for each $7\leq l \leq N$:
\begin{equation}\label{J-function2}
	\left.\begin{array}{l}
		-\omega_2J_l=H_{l+2}+\partial_J H_1J_{l-2}+\sum\limits_{k=5}^{l}\partial_J H_kJ_{l-k+2}+\dfrac{\beta_2}{2}J^2_{\lfloor\tfrac{l}{2}\rfloor-1}\delta_{\lfloor\tfrac{l}{2}\rfloor,\tfrac{l}{2}}\\ \qquad\qquad +\beta_2\sum\limits_{k=3}^{\lfloor\tfrac{l-1}{2}\rfloor}J_{k-1}J_{l-k-1}+\sum\limits_{k=5}^{l-2}\sum\limits_{\gamma\in\widetilde{\Gamma}_{l-k+2}}\dfrac{\partial^{|\gamma|}_JH_{k}}{\gamma!}\prod\limits_{j=1}^{l-k-1}J_{j+1}^{\gamma_{j}},
	\end{array}\right.
\end{equation}
where $\delta_{j,k}$ denotes the Kronecker's Delta Function.

	Observe that the functions in \eqref{J-function1}-\eqref{J-function2} are polynomials with respect to the variables \eqref{levelSet-AngularVariables}, all of which are defined recursively, while $J_2$ has no angular dependence. Thus, when $I$ is fixed, $J_2$ becomes constant and $\{J_l:\ l\geq3\}$ is a set of polynomials in the angular variables defined in \eqref{levelSet-AngularVariables}. In the following lemma we compute the degree of the polynomials $J_l$, which we denote by ${\dd}[f]$.

\begin{Lem}\label{J-Lemma}
	For any $l\geq4$, it holds that
	\begin{equation}\label{InductiveOrder}
		\left\{\begin{array}{ll}
			{\dd}[J_l]={\dd}[J_{l-1}]+3,&\text{if}\quad [l]_{\ZZ_3}=0.\\
			{\dd}[J_l]={\dd}[J_{l-1}]+1, & \text{otherwise}.
		\end{array}\right.
	\end{equation}
	In particular,
	\begin{equation}\label{J-Degree}
		{\dd}[J_l]=l+2\left\lfloor\tfrac{l}{3}\right\rfloor,\quad l\geq3.
	\end{equation}
\end{Lem}
\proof
It is immediate to check that
\begin{equation}
	\left.\begin{array}{ccccc}
		\dd[J_2]=0,&\dd[J_3]=5,&\dd[J_4]=6,&\dd[J_5]=7,& \dd[J_6]=10.
	\end{array}\right.
\end{equation}
Moreover, the next identity holds for  $l\in\{7,8,9\}$:
\begin{equation}\label{DegreeExample}
	\dd[J_l]=\dd\left[\sum_{k\geq5}^{l}\partial_J H_kJ_{l-k+2}+\sum\limits_{k=5}^{l-2}\sum\limits_{\gamma\in\widetilde{\Gamma}_{l-k+2}}\dfrac{\partial^{|\gamma|}_JH_{k}}{\gamma!}\prod_{j=1}^{l-k-1}J_{j+1}^{\gamma_{j}}\right].
\end{equation} 
Concerning this,
\begin{equation}
	\dd[J_l]>\max\limits_{\gamma\in\widetilde{\Gamma}_l}\ \dd\left[\prod_{j=1}^{l-3}J_{j+1}^{\gamma_j}\right],\ \hbox{when }l\in\{4,5\},
\end{equation}which implies
\begin{equation}
	\dd[J_7]=\dd\left[\sum\limits_{k=5}^6\partial_J H_kJ_{9-k}\right]=11,\ \hbox{ and }\
	\dd[J_8]=\dd\left[\sum\limits_{k=5}^7\partial_J H_kJ_{10-k}\right]=12,
\end{equation}
where we have not taken into account the terms from  \eqref{DegreeExample} with non maximum degree. On the other hand,
\begin{equation}
	\dd[J_6]=\max\limits_{\gamma\in\widetilde{\Gamma}_6}\ \dd\left[\prod_{j=1}^{3}J_{j+1}^{\gamma_j}\right]=\dd\left[J_{3}^{2}\right]=2\dd\left[J_3\right]=10,
\end{equation} and
\begin{equation}
	\dd[J_9]=\dd\left[\partial_J H_5J_{6}+\sum\limits_{\gamma\in\widetilde{\Gamma}_{6}}\dfrac{\partial^{|\gamma|}_JH_{5}}{\gamma!}\prod_{j=1}^{3}J_{j+1}^{\gamma_{j}}\right]=\dd\left[\partial_J H_5J_{6}+\tfrac{1}{2}\partial^{2}_JH_{5}J_3^2\right]=15.
\end{equation}
Proceeding analogously for $l\in\{10,11\}$, the claim of the lemma is proved for $4\leq l\leq11$:
\begin{equation}\label{Examples}
	\left.\begin{array}{lll}
		\dd[J_3]=5,&\dd[J_4]=6,&\dd[J_5]=7,\\ \dd[J_6]=10,& \dd[J_7]=11,&\dd[J_8]=12,\\ \dd[J_9]=15,&\dd[J_{10}]=16,& \dd[J_{11}]=17.
	\end{array}\right.
\end{equation}
Moreover, observe that we have established the following inequality for each $7\leq l\leq11$:
\begin{equation}\label{J-property}
	\dd[J_{k}]\geq\dd\left[\prod_{j=1}^{k-3}J_{j+1}^{\gamma_j}\right],\quad\hbox{ for all }\gamma\in\widetilde{\Gamma}_{k},\ \hbox{ with }4\leq k\leq l-3,
\end{equation}
being strict the inequality when $k\in\{4,5\}$. Next, we shall use this inequality to establish~\eqref{InductiveOrder} by induction.

Let $l\in\mathbb{N}$ be large enough such that \eqref{InductiveOrder} is satisfied for all $3\leq k\leq l-1$, and \eqref{J-property} holds for all $6\leq k\leq l-4$. First, it is immediate to obtain \eqref{DegreeExample} for the case $k=l$, i.e.,
\begin{equation}
	\dd[J_l]=\dd\left[\sum_{k\geq5}^{l}\partial_J H_kJ_{l-k+2}+\sum\limits_{k=5}^{l-2}\sum\limits_{\gamma\in\widetilde{\Gamma}_{l-k+2}}\dfrac{\partial^{|\gamma|}_JH_{k}}{\gamma!}\prod_{j=1}^{l-k-1}J_{j+1}^{\gamma_{j}}\right].
\end{equation} 
Then, assuming that \eqref{J-property} holds for $k=l-3$, we obtain
\begin{equation}
	\dd[J_l]=\dd\left[\sum_{k\geq5}^{5+[l]_{\ZZ_3}}\partial_J H_kJ_{l-k+2}\right]=5+\dd[J_{l-3}],
\end{equation} 
which, after considering the different cases according to the value of $[l]_{\ZZ_3}$, yields the desired formula~\eqref{InductiveOrder} for the general case $k=l$. To complete the argument, it then remains to prove \eqref{J-property} for the case $k=l-3$. 

To this end, for any $n\leq\lfloor\tfrac{k-3}{2}\rfloor$, let $\gamma^n\in\mathbb{N}^{k-3}$ be the vector of components
\begin{equation}\label{vector-proof}
	\gamma_j=\delta_{jm}, \quad\hbox{where }\quad m\in\{n,k-2-n\}.
\end{equation}
Moreover, when $k$ is even, we also consider the case $n=\tfrac{k-2}{2}$, defined as $$\gamma_{j}=2\delta_{jn},\qquad1\leq j\leq k-3.$$
Observe that $\gamma^n\in\widetilde{\Gamma}_k$, and
\begin{equation}
	\prod_{j=1}^{k-3}J_{j+1}^{\gamma_{j}}=J_{n+1}J_{k-n-1},\quad\hbox{for all}\ n\leq\lfloor\tfrac{k-3}{2}\rfloor.
\end{equation}
Thus, by  applying standard inequalities for the floor function in \eqref{J-Degree}, we obtain
\begin{align}\label{J-Inequality}
	\dd\left[J_{n+1}J_{k-n-1}\right]=\dd[J_{n+1}]+\dd\left[J_{k-n-1}\right]\leq \dd[J_k],\ \hbox{ for all}\ n\leq\lfloor\tfrac{k-3}{2}\rfloor,
\end{align}
which also holds for the case $n=\tfrac{k-2}{2}$ with $k$ even. On the other hand, if $\gamma\in\widetilde{\Gamma}_k$ is different from these vectors $\gamma^n$, let $\mathbf{n}$ be the lowest sub-index $j$ for which $\gamma_{j}\neq0$. In that case, 
\begin{equation}
	k=\sum\limits_{j}\gamma_j (j+1)=\mathbf{n}+1+\sum\limits_{j}\widetilde{\gamma}_j (j+1),\quad\hbox{where }\widetilde{\gamma}_j=\gamma_j-\delta_{j\mathbf{n}}.
\end{equation}
Notice that necessarily $\widetilde{\gamma}_j=0$, for any $j>k-\mathbf{n}-3$, so it can be asserted without any loss of generality that $\widetilde{\gamma}\in\widetilde{\Gamma}_{k-\mathbf{n}-1}$. Then, we conclude
\begin{align}
	\dd\left[\prod_{j=1}^{k-3}J_{j+1}^{\gamma_j}\right]&=\sum\limits_{j}\gamma^j\dd[J_{j+1}]=\dd[J_{\mathbf{n}+1}]+\dd\left[\prod_{j=1}^{k-\mathbf{n}-4}J_{j+1}^{\widetilde{\gamma}_j}\right]\\ &\leq\dd[J_{\mathbf{n}+1}]+\dd[J_{k-\mathbf{n}-1}]\leq \dd[J_{k}],
\end{align}
which completes the proof.
\cqfd 

The previous lemma has a corollary that is crucial in order to prove Proposition~\ref{J-Proposition}:

\begin{Cor}\label{J-Corollary}
	For any $l\geq4$, let $\widetilde{\Gamma}^M_{l}\subset\widetilde{\Gamma}_{l}$ be the subset of elements with maximum degree, i.e.,
	\begin{equation}\label{J-MaxSet}
		\widetilde{\Gamma}^M_{l}=\left\{\gamma\in\widetilde{\Gamma}_{l}:{\dd}\left[\prod_{k}J_{k+1}^{\gamma_k}\right]={\dd}[J_l]\right\}.
	\end{equation}
	Then, $\widetilde{\Gamma}^M_{4}=\widetilde{\Gamma}^M_{5}=\emptyset$, while $\widetilde{\Gamma}^M_{l}$ has the following structure for $l\geq6$:
	\begin{itemize}
		\item [i)] Case $[l]_{\ZZ_3}=0$:
		\begin{equation}
			\gamma\in\widetilde{\Gamma}^M_{l} \ \Leftrightarrow\ 
			\gamma_j=0,\text{ if }j\neq 3k-1,\ k\geq1.
		\end{equation}
		\item [ii)] Case $[l]_{\ZZ_3}=1$:
		\begin{equation}
			\gamma\in\widetilde{\Gamma}^M_{l} \ \Leftrightarrow\ \gamma:\left\{\begin{array}{cc}
				\gamma_{\tilde{j}}=1,&\text{for a unique }\tilde{j}=3k,\ k\geq1,\\
				\gamma_j=0,&\text{if }j\neq 3k-1,\ k\geq1, \hbox{ and }j\neq \tilde{j}.
			\end{array}\right.
		\end{equation}
		\item [iii)] When $[l]_{\ZZ_3}=2$, each vector of $\widetilde{\Gamma}^M_{l}$ has one of the following structures. Either $\gamma\in\widetilde{\Gamma}^M_{l}$  is such that
		\begin{equation}
			\gamma:\left\{\begin{array}{cc}
				\gamma_{\tilde{j}}=1,&\text{for a unique }\tilde{j}=3k+1,\ k\geq1,\\
				\gamma_j=0,&\text{if }j\neq 3k-1,\ k\geq1, \hbox{ and }j\neq \tilde{j},
			\end{array}\right.
		\end{equation}
		or there exist indices $\{j_1,j_2\}$, not necessarily different, such that $[j_1]_{\ZZ_3}=[j_2]_{\ZZ_3}=1$, and
		\begin{equation}
			\gamma:\left\{\begin{array}{lc}
				\gamma_{j}=1+\delta_{j_1j_2},&\text{if }j\in\{j_1,j_2\},\\
				\gamma_j=0,&\text{if }j\neq 3k-1,\ k\geq1, \hbox{ and }j\notin\{j_1,j_2\}.
			\end{array}\right.
		\end{equation}
	\end{itemize}
\end{Cor}
\proof
Since the assertion $\widetilde{\Gamma}^M_{4}=\widetilde{\Gamma}^M_{5}=\emptyset$ is trivial, from now on we assume $l\geq6$. We also recall that any natural number admits the decomposition $l=3\left\lfloor\tfrac{l}{3}\right\rfloor+ [l]_{\ZZ_3},$ so \eqref{J-Degree} can be written as
\begin{equation}
	\dd[J_l]=l+\dfrac{2}{3}\left(l-[l]_{\ZZ_3}\right),\ \hbox{ for any }l\geq3.
\end{equation}
Moreover, we consider a vector $\gamma^j\in\widetilde{\Gamma}_l$, with $j>1$, defined as in \eqref{vector-proof}. Then,
\begin{equation}
	\dd[J_{j+1}]+\dd[J_{l-j-1}]=l+\dfrac{2}{3}\left(l-[j+1]_{\ZZ_3}-[l-j-1]_{\ZZ_3}\right),
\end{equation}
and
\begin{equation}\label{J-Equality}
	\gamma^j\in\widetilde{\Gamma}_l^M\Leftrightarrow [l]_{\ZZ_3}=[j+1]_{\ZZ_3}+[l-j-1]_{\ZZ_3}.
\end{equation} 
Therefore, when $[l]_{\ZZ_3}=0$, it follows directly that
\begin{align}
	\gamma^j\in\widetilde{\Gamma}_l^M\Leftrightarrow \ [j+1]_{\ZZ_3}=0.
\end{align}
Using this, the implication to the left in $i)$ is obtained by recurrence, while the other can be easily proved by contradiction. In particular,
\begin{align}
	\dd[J_l]=\dd\left[J_3^{\left\lfloor\tfrac{l}{3}\right\rfloor}\right]=\left\lfloor\tfrac{l}{3}\right\rfloor\dd\left[J_3\right],\ \hbox{when }[l]_{\ZZ_3}=0.
\end{align}
In addition, the procedure is analogous for the case $ii)$, that is $[l]_{\ZZ_3}=1$, which leads to
\begin{align}
	\dd[J_l]=\dd\left[J_3^{\left\lfloor\tfrac{l}{3}\right\rfloor-1}J_4\right]=\left(\left\lfloor\tfrac{l}{3}\right\rfloor-1\right)\dd\left[J_3\right]+\dd\left[J_4\right],\ \hbox{when }[l]_{\ZZ_3}=1.
\end{align}
On the other hand,  according to \eqref{J-Equality}, in the third case there are two kinds of vector in order to obtain $\gamma^j\in\widetilde{\Gamma}_l^M$. More concretely, these are the following:
\begin{align}
	[j+1]_{\ZZ_3}=0\Rightarrow[l-j-1]_{\ZZ_3}=[l]_{\ZZ_3},\ \hbox{ and }\ [j+1]_{\ZZ_3}=1\Rightarrow[l-j-1]_{\ZZ_3}=1.
\end{align}
Furthermore, when $l\geq8$ and $[l]_{\ZZ_3}=2$, it is also satisfied that
\begin{align}
	\dd[J_l]=\dd\left[J_3^{\left\lfloor\tfrac{l}{3}\right\rfloor-1}J_5\right]=\left(\left\lfloor\tfrac{l}{3}\right\rfloor-1\right)\dd\left[J_3\right]+\dd\left[J_5\right],
\end{align}
and
\begin{align}
	\dd[J_l]=\dd\left[J_3^{\left\lfloor\tfrac{l}{3}\right\rfloor-2}J_4^2\right]=\left(\left\lfloor\tfrac{l}{3}\right\rfloor-2\right)\dd\left[J_3\right]+2\dd\left[J_4\right].
\end{align}
From this, and using the previous cases, $iii)$ is easily obtained after proceeding again by recurrence and contradiction.
\cqfd

\subsection{Proof of Proposition~\ref{J-Proposition}}
For any $l\geq7$, by \eqref{DegreeExample} and \eqref{J-MaxSet} we can write
\begin{equation}
	J_l=-\dfrac{1}{\omega_2}\left(\sum_{k=5}^{5+[l]_{\ZZ_3}}\partial_J H_kJ_{l-k+2}+\sum\limits_{k=5}^{5+[l]_{\ZZ_3}}\sum\limits_{\gamma\in\widetilde{\Gamma}^M_{l-k+2}}\dfrac{\partial^{|\gamma|}_JH_{k}}{\gamma!}\prod_{j=1}^{l-k-1}J_{j+1}^{\gamma_{j}}\right)+\text{l.o.t.}
\end{equation}
In particular, when $[l]_{\ZZ_3}=0$ this implies that
\begin{equation}\label{L0H5}
	J_l=-\dfrac{1}{\omega_2}\left(\partial_J H_5 J_{l-3}+\sum\limits_{\gamma\in\widetilde{\Gamma}^M_{l-3}}\dfrac{\partial^{|\gamma|}_JH_{5}}{\gamma!}\prod_{j=1}^{l-6}J_{j+1}^{\gamma_{j}}\right)+\text{l.o.t.}
\end{equation}
and, since $[l-3]_{\ZZ_3}=[l]_{\ZZ_3}$, we can iterate using the recursion with $J_{l-3}$ to obtain
\begin{align}
	J_l=&\ \dfrac{1}{\omega^2_2}\left((\partial_J H_5)^2 J_{l-6}+\partial_J H_5\sum\limits_{\gamma\in\widetilde{\Gamma}^M_{l-6}}\dfrac{\partial^{|\gamma|}_JH_{5}}{\gamma!}\prod_{j=1}^{l-9}J_{j+1}^{\gamma_{j}}\right)\\&-\dfrac{1}{\omega_2}\sum\limits_{\gamma\in\widetilde{\Gamma}^M_{l-3}}\dfrac{\partial^{|\gamma|}_JH_{5}}{\gamma!}\prod_{j=1}^{l-6}J_{j+1}^{\gamma_{j}}+\text{l.o.t.}
\end{align}
This understood, \eqref{L0H5} goes into the following expression:
\begin{align}
	J_l=&\left(\dfrac{-1}{\omega_2}\right)^{\tfrac{l}{3}-1}\left(\partial_J H_5\right)^{\tfrac{l}{3}-1} J_{3}-\sum_{k=1}^{\tfrac{l}{3}-2}\dfrac{\left(-\partial_JH_5\right)^{k-1}}{\omega_2^k}\sum\limits_{\gamma\in\widetilde{\Gamma}^M_{l-3k}}\dfrac{\partial^{|\gamma|}_JH_{5}}{\gamma!}\prod_{j=1}^{l-3(k+1)}J_{j+1}^{\gamma_{j}}+\text{l.o.t.}\\
	=&\left(\dfrac{-1}{\omega_2}\right)^{\tfrac{l}{3}}\left(\partial_J H_5\right)^{\tfrac{l}{3}-1}H_5-\sum_{k=1}^{\tfrac{l}{3}-2}\dfrac{\left(-\partial_JH_5\right)^{k-1}}{\omega_2^k}\sum\limits_{\gamma\in\widetilde{\Gamma}^M_{l-3k}}\dfrac{\partial^{|\gamma|}_JH_{5}}{\gamma!}\prod_{j=1}^{l-3(k+1)}J_{j+1}^{\gamma_{j}}+\text{l.o.t.}
\end{align}
Clearly, the first term belongs to $\mathcal{L}^5_0\left[\tfrac{l}{3},1\right]$. On the other hand, proceeding by induction and using $i)$ from Corollary~\ref{J-Corollary}, it follows that
\begin{equation}
	\prod_{j=1}^{l-3(k+1)}J_{j+1}^{\gamma_{j}}=\left(\dfrac{-1}{\omega_2}\right)^{\sum\limits_j \tfrac{\gamma_j(j+1)}{3}}\prod\limits_{j}\mathcal{P}_j^{\gamma_{j}}+\hbox{l.o.t.}=\left(\dfrac{-1}{\omega_2}\right)^{ \tfrac{l}{3}-k}\prod\limits_{j}\mathcal{P}_j^{\gamma_{j}}+\hbox{l.o.t.}
\end{equation} 
for all $1\leq k\leq\tfrac{l}{3}-2$, and all $\gamma\in\widetilde{\Gamma}^M_{l-3k}$. Moreover,
\begin{equation*}
	\mathcal{P}^{\gamma_{j}}_j\in\mathcal{L}_0^5\left[\gamma_{j}\tfrac{j+1}{3},\tfrac{j+1}{3}-1\right],\quad\hbox{for all }1\leq j\leq l-3(k+1).
\end{equation*} Then, using the inclusion 
\begin{equation*}
\mathcal{L}_0^5[\cdot,n]\subset\mathcal{L}_0^5[\cdot,m],\quad\hbox{when }n<m,
\end{equation*}
we conclude that
\begin{equation}
	\prod\limits_{j}\mathcal{P}_j^{\gamma_{j}}\in\mathcal{L}_0^5\left[\sum\limits_j \tfrac{\gamma_j(j+1)}{3},\max\limits_{\gamma\in\widetilde{\Gamma}^M_{l-3}}|\gamma|\right]=\mathcal{L}_0^5\left[\tfrac{l}{3}-k,\tfrac{l}{3}-1\right],
\end{equation}
for all $1\leq k\leq\tfrac{l}{3}-2$, and all $\gamma\in\widetilde{\Gamma}^M_{l-3k}$. Therefore, for every $1\leq k\leq\tfrac{l}{3}-2$, the higher order terms of 
\begin{equation}
	\left(\partial_JH_5\right)^{k-1}\sum\limits_{\gamma\in\widetilde{\Gamma}^M_{l-3k}}\dfrac{\partial^{|\gamma|}_JH_{5}}{\gamma!}\prod_{j=1}^{l-3(k+1)}J_{j+1}^{\gamma_{j}},
\end{equation}
belong to $\mathcal{L}_0^5\left[\tfrac{l}{3},\tfrac{l}{3}-1\right]$, which proves $i)$. When $[l]_{\ZZ_3}=1$, we have
\begin{align}
	J_l=&-\dfrac{1}{\omega_2}\left(\partial_J H_5 J_{l-3}+\sum\limits_{\gamma\in\widetilde{\Gamma}^M_{l-3}}\dfrac{\partial^{|\gamma|}_JH_{5}}{\gamma!}\prod_{j=1}^{l-6}J_{j+1}^{\gamma_{j}}\right)\\
	&-\dfrac{1}{\omega_2}\left(\partial_J H_6 J_{l-4}+\sum\limits_{\gamma\in\widetilde{\Gamma}^M_{l-4}}\dfrac{\partial^{|\gamma|}_JH_{6}}{\gamma!}\prod_{j=1}^{l-7}J_{j+1}^{\gamma_{j}}\right)+\text{l.o.t.}
\end{align}
Observe that the second line satisfies $ii)$ because $[l-4]_{\ZZ_3}=0$, thus the case $i)$ applies directly to $J_{l-4}$ and $\prod\limits_{j}J_{j+1}^{\gamma_{j}}$, for all $\gamma\in\widetilde{\Gamma}^M_{l-4}$. Concerning the remainders, since $[l-3k]_{\ZZ_3}=[l]_{\ZZ_3}=1$, for all $k\in\mathbb{N}$, we can iterate with every function $J_{l-3k}$, and hence the first line is rewritten as follows:
\begin{align}
	\left(\dfrac{-1}{\omega_2}\right)^{\tfrac{l-1}{3}}\partial_J H_5^{\tfrac{l-1}{3}-1}H_6-\sum_{k=1}^{\tfrac{l-1}{3}-2}\left(\dfrac{-\partial_JH_5}{\omega_2}\right)^{k-1}\sum\limits_{\gamma\in\widetilde{\Gamma}^M_{l-3k}}\dfrac{\partial^{|\gamma|}_JH_{5}}{\gamma!}\prod_{j=1}^{l-3(k+1)}J_{j+1}^{\gamma_{j}}.
\end{align}
Finally, for all $\gamma\in\widetilde{\Gamma}^M_{l-3k}$, with $1\leq k\leq \tfrac{l-1}{3}-2$, its decomposition according to assertion $ii)$ in Corollary~\ref{J-Corollary} implies that 
\begin{align}
	\prod_{j=1}^{l-3(k+1)}J_{j+1}^{\gamma_{j}}=\left(\dfrac{-1}{\omega_2}\right)^{\tfrac{l-1}{3}-k}\sum_{n=0}^{\left\lfloor\tfrac{l}{3}\right\rfloor-(k+2)}\partial_J^nH_6\mathcal{P}_n+\hbox{ l.o.t.}\ 
\end{align}
where $\mathcal{P}_n\in\mathcal{L}_0^5\left[\tfrac{l-1}{3}-(k+1),\tfrac{l-1}{3}-k\right],$ for all $n$, and proceeding as above, $ii)$ is obtained. To conclude, when $[l]_{\ZZ_3}=2$ we have
\begin{align}
	J_l=&-\dfrac{1}{\omega_2}\left(\partial_J H_5 J_{l-3}+\sum\limits_{\gamma\in\widetilde{\Gamma}^M_{l-3}}\dfrac{\partial^{|\gamma|}_JH_{5}}{\gamma!}\prod_{j=1}^{l-6}J_{j+1}^{\gamma_{j}}\right)\\
	&-\dfrac{1}{\omega_2}\left(\partial_J H_6 J_{l-4}+\sum\limits_{\gamma\in\widetilde{\Gamma}^M_{l-4}}\dfrac{\partial^{|\gamma|}_JH_{6}}{\gamma!}\prod_{j=1}^{l-7}J_{j+1}^{\gamma_{j}}\right)\\
	&-\dfrac{1}{\omega_2}\left(\partial_J H_7 J_{l-5}+\sum\limits_{\gamma\in\widetilde{\Gamma}^M_{l-5}}\dfrac{\partial^{|\gamma|}_JH_{7}}{\gamma!}\prod_{j=1}^{l-8}J_{j+1}^{\gamma_{j}}\right)+\text{l.o.t.}
\end{align}
In this case, $[l-5]_{\ZZ_3}=0$ and $[l-4]_{\ZZ_3}=1$. Then, by applying the cases $i)$ and $ii)$, it follows that the sum of the higher order terms of the last two lines belongs to $\mathcal{L}_2^7\left[\left\lfloor\tfrac{l}{3}\right\rfloor,\left\lfloor\tfrac{l}{3}\right\rfloor-1\right]$. Using this, since $[l-3k]_{\ZZ_3}=[l]_{\ZZ_3}=2$, for all $k\in\mathbb{N}$, we can iterate with each $J_{l-3k}$, and because of the structure of $\widetilde{\Gamma}^M_{l-3k}$ (see $iii)$ in Corollary~\ref{J-Corollary}), it follows that the higher order term in the first line is an element of $\mathcal{L}_2^7\left[\left\lfloor\tfrac{l}{3}\right\rfloor,\left\lfloor\tfrac{l}{3}\right\rfloor-1\right]$, and we are done.

\section{Proof of Proposition~\ref{Pro-Hyp0}}\label{Section-Proof1}
Our first goal is to provide the explicit dependence with respect to $\varepsilon$ of the numerator and the inverse of the denominator in
\begin{equation*}
	\sum\limits_{l=3}^{N}\varepsilon^{l}\dfrac{\partial_\varphi H_{l+2}(I,J,\vec{\varphi})}{\partial_J H(I,J,\vec{\varphi})},
\end{equation*}
presenting each of them in Propositions~\ref{Pro-Hl-numerator} and~\ref{Pro-Hl-denominator}, respectively. Proposition~\ref{Pro-Hyp0} will then follow as an easy corollary.
\begin{Pro}\label{Pro-Hl-numerator}
	Let $V\in\mathcal{V}_0$ be satisfying \eqref{nonresonant}. Then, on every level set \eqref{levelSet} with sufficiently small energy $h>0$, it holds that
	\begin{equation*}
		\sum\limits_{l=3}^{N}\varepsilon^{l}\partial_\varphi H_{l+2}(I,J(I,\bar{\varphi},\varepsilon),\bar{\varphi})=\sum\limits_{l=3}^{N}\varepsilon^l h^\varphi_l(I,J_0(I),\vec{\varphi})+\Theta(\varepsilon^{N+1}),
	\end{equation*}
	{with arbitrary $N\geq3$, and} where each $h^\varphi_l$ is a polynomial of degree $l+2\left\lfloor\tfrac{l}{3}\right\rfloor$ in the angular variables defined in \eqref{levelSet-AngularVariables}. More concretely,
	\begin{equation}
		h_l^\varphi=\left(\dfrac{-1}{\omega_2}\right)^{\left\lfloor\tfrac{l}{3}\right\rfloor-1}\hh^\varphi_l\ +\hbox{ l.o.t.}
	\end{equation}
	for some explicit $\hh^\varphi_l\in\mathcal{L}_{[l]_{\ZZ_3},\varphi}^{5+[l]_{\ZZ_3}}\left[\left\lfloor\tfrac{l}{3}\right\rfloor,\left\lfloor\tfrac{l}{3}\right\rfloor-1\right].$ 
\end{Pro}
\proof
{Let us assume $N\geq7$, otherwise the proof is trivial}. From Section~\ref{Section-ReversedAction} we know that
\begin{align}
	\sum_{l\geq3}^N\varepsilon^{l}\partial_\varphi H_{l+2}(I,J,\vec{\varphi})=&\ \sum_{l\geq3}^N\varepsilon^{l}\partial_\varphi H_{l+2}+\sum_{l\geq3}^{N-2}\sum\limits_{k=2}^{N-l}\varepsilon^{l+k}\partial^{J}_\varphi H_{l+2}J_{k}\\&\ +\sum_{l\geq3}^{N-4}\sum\limits_{k=4}^{N-l}\varepsilon^{l+k}\sum\limits_{\gamma\in\widetilde{\Gamma}_k}\dfrac{\partial^{J,|\gamma|}_\varphi H_{l+2}}{\gamma!}\prod_{j=1}^{k-3}J_{j+1}^{\gamma_{j}}+\Theta(\varepsilon^{N+1}),
\end{align}
with all terms evaluated at $(I,J_0,\vec{\varphi})$. We now distinguish three different cases:\vspace{0.1cm}
\begin{itemize}
	\item For $l\in\{3,4\}$, $h^\varphi_l=\partial_\varphi H_{l+2}$, while $h^\varphi_5=\partial_\varphi H_{7}+\hbox{l.o.t.}$\vspace{0.3cm}
	\item For $l\in\{6,7,8\}$:
	\begin{equation*}
		h^\varphi_l=\sum_{k=0}^{[l]_{\ZZ_3}}\partial^J_\varphi H_{5+k} J_{l-3-k}+\hbox{l.o.t.}
	\end{equation*}
	\item For $l\geq9$:
	\begin{equation*}
		h^\varphi_l=\sum_{k=0}^{[l]_{\ZZ_3}}\left[\partial^J_{\varphi} H_{5+k} J_{l-3-k}+\sum\limits_{\gamma\in\widetilde{\Gamma}^M_{l-3-k}}\dfrac{\partial^{J,|\gamma|}_{\varphi}H_{5+k}}{\gamma!}\prod_{j=1}^{l-6-k}J_{j+1}^{\gamma_{j}}\right]+\hbox{l.o.t.}
	\end{equation*}
	
\end{itemize}
Then, by analyzing each term in these expressions, the result follows directly from Proposition~\ref{J-Proposition}.
\cqfd
\begin{Pro}\label{Pro-Hl-denominator}
	Let $V\in\mathcal{V}_0$ be satisfying \eqref{nonresonant}. Then, on every level set \eqref{levelSet} with sufficiently small energy $h>0$, it holds that
	\begin{align}
		\dfrac{1}{\partial_J H(I,J,\vec{\varphi})}=&\dfrac{1}{\omega_2}-\varepsilon^2\dfrac{\partial_J H_1(I,J_0)}{\omega_2^2}+\sum_{l=3}^{N}\varepsilon^lh_l(I,J_0,\vec{\varphi})+\Theta(\varepsilon^{N+1}),
	\end{align}
	{with arbitrary $N\geq3$, and} where each $h_l$ is a polynomial of degree $l+2\left\lfloor\tfrac{l}{3}\right\rfloor$ in the angular variables defined in \eqref{levelSet-AngularVariables}. More concretely,
	\begin{equation*}
		h_l=-\left(\dfrac{-1}{\omega_2}\right)^{\left\lfloor\tfrac{l}{3}\right\rfloor+1}\hh_l+\textnormal{ l.o.t.}
	\end{equation*}
	for some explicit $\hh_l\in\mathcal{L}_{[l]_{\ZZ_3}}^{5+[l]_{\ZZ_3}}\left[\left\lfloor\tfrac{l}{3}\right\rfloor,1\right].$
\end{Pro}
\proof
We begin by recalling that $h_3=-\omega_2^{-2}\partial_J H_5$, and we claim that
\begin{align}\label{hl-formula}
	h_l=\dfrac{-1}{\omega_2^2}\partial_JH_{l+2}+\sum_{\gamma\in\widetilde{\Gamma}_l}(-1)^{|\gamma|}\dfrac{|\gamma|!}{\gamma!}\dfrac{1}{\omega_2^{|\gamma|+1}}\left(2\partial_J H_1\right)^{\gamma_1}\prod_{k=2}^{l-3}\left(\partial_{J}H_{k+3}\right)^{\gamma_k},
\end{align} for $l\geq4$, which will be proved at the end.

First, for $l\in\{3,4,5\}$, the higher order term of $h_l$ is uniquely determined by $\partial_J H_{l+2}$. Nevertheless, when $l\geq6$, the product at \eqref{hl-formula} implies that
\begin{equation}\label{degree}
	{\dd}[h_l]=\sum_{k=2}^{l-3}\gamma_k(k+3)=l+2\max\limits_{\gamma\in\widetilde{\Gamma}_l}\sum_{k=2}^{l-3}\gamma_k.
\end{equation}
Therefore, we have to check the subset of $\widetilde{\Gamma}_l$ which maximizes this quantity. To this end, when $[l]_{\ZZ_3}=0$, it is immediate to see that there exists a unique $\gamma\in\widetilde{\Gamma}_l$ maximizing \eqref{degree}, defined by
\begin{equation*}
	\gamma_2=\left\lfloor\frac{l}{3}\right\rfloor,\quad\hbox{and}\quad \gamma_j=0,\hbox{ for }j\neq2.
\end{equation*} Similarly, for $[l]_{\ZZ_3}=1$ we obtain again uniqueness of $\gamma$, which is given by 
\begin{equation*}
	\gamma_2=\left\lfloor\frac{l}{3}\right\rfloor-1,\quad\gamma_3=1,\quad\hbox{and}\quad\gamma_j=0,\quad\hbox{ for }j\notin\{2,3\}.
\end{equation*}However, if $[l]_{\ZZ_3}=2$, then there are two vectors in $\widetilde{\Gamma}_l$ maximizing \eqref{degree}. These are defined by
\begin{equation*}
	\gamma\in\widetilde{\Gamma}_l:\quad\gamma_2=\left\lfloor\tfrac{l}{3}\right\rfloor-2,\quad\gamma_3=1,\quad\hbox{and}\quad \gamma_j=0,\quad\hbox{for } j\notin\{2,3\},
\end{equation*}
and
\begin{equation*}
	\gamma\in\widetilde{\Gamma}_l:\quad\gamma_2=\left\lfloor\tfrac{l}{3}\right\rfloor-1,\quad\gamma_4=1,\quad\hbox{and}\quad \gamma_j=0,\quad\hbox{for } j\notin\{2,4\}.
\end{equation*}
Consequently, using the expressions of these vectors in \eqref{hl-formula} gives the following structure for $h_l$, which satisfies the statement of the proposition:
\begin{itemize}
	\item [i)] Case $[l]_{\ZZ_3}=0$:
	$$h_l=-\left(\dfrac{-1}{\omega_2}\right)^{\left\lfloor\tfrac{l}{3}\right\rfloor+1}\left(\partial_{J}H_{5}\right)^{\left\lfloor\tfrac{l}{3}\right\rfloor}+\textnormal{ l.o.t.}$$
	\item [ii)]Case $[l]_{\ZZ_3}=1$:
	$$h_l=-\left\lfloor\dfrac{l}{3}\right\rfloor\left(\dfrac{-1}{\omega_2}\right)^{\left\lfloor\tfrac{l}{3}\right\rfloor+1}\left(\partial_{J}H_{5}\right)^{\left\lfloor\tfrac{l}{3}\right\rfloor-1}\partial_{J}H_{6}+\textnormal{ l.o.t.}$$
	\item [ii)]Case $[l]_{\ZZ_3}=2$:
	$$h_l=-\left\lfloor\dfrac{l}{3}\right\rfloor\left(\dfrac{-1}{\omega_2}\right)^{\left\lfloor\tfrac{l}{3}\right\rfloor+1}\left(\partial_{J}H_{5}\right)^{\left\lfloor\tfrac{l}{3}\right\rfloor-2}\left[\partial_{J}H_{5}\partial_{J}H_{7}+\dfrac{\left\lfloor\tfrac{l}{3}\right\rfloor-1}{2}\left(\partial_{J}H_{6}\right)^2\right]+\textnormal{ l.o.t.}$$
\end{itemize}
To conclude, it remains to prove Equation~\eqref{hl-formula}, for which we recall the following identities from the proof of Proposition~\ref{Pro-Continuous}:
\begin{align}
	\left.\dfrac{d^2}{d\varepsilon^2}\partial_{J}H\right|_{\varepsilon=0}=2\partial_J H_1(I,J_0),\qquad\left.\dfrac{d^3}{d\varepsilon^3}\partial_{J}H\right|_{\varepsilon=0}=6\partial_J H_5(I,J_0,\vec{\varphi}).
\end{align}
In addition, the generalized chain rule formula gives that
\begin{align}
	\left.\dfrac{d^l}{d\varepsilon^l}\partial_{J}H\right|_{\varepsilon=0}=&\left.\left(\partial^l_\varepsilon \partial_JH+\partial^{2}_JH\partial_l J+\sum_{\gamma\in\widetilde{\Gamma}_l}\dfrac{l!}{\gamma!}\partial_J^{|\gamma|+1}H\prod_{k=1}^{l-3}\left(\dfrac{\partial_{k+1} J}{(k+1)!}\right)^{\gamma_{k}}\right)\right|_{\varepsilon=0}\\&=l!\partial_JH_{l+2}(I,J_0,\vec{\varphi}), \quad \hbox{ for }l\geq4.
\end{align}
On the other hand, the lower order terms in the Taylor expansion of $(\partial_J H)^{-1}$ at $\varepsilon=0$ are determined by:
\begin{align}
	\left.\dfrac{d}{d\varepsilon}\dfrac{1}{\partial_J H}\right|_{\varepsilon=0}=\left.\dfrac{-1}{(\partial_J H)^2}\left(\partial_\varepsilon \partial_J H+\partial^2_J H \partial_\varepsilon L\right)\right|_{\varepsilon=0}=0,
\end{align}
and
\begin{align}
	\left.\dfrac{d^2}{d\varepsilon^2}\dfrac{1}{\partial_J H}\right|_{\varepsilon=0}=&\dfrac{-\partial_J H_1(I,J_0)}{\omega_2^2},\qquad\left.\dfrac{d^3}{d\varepsilon^3}\dfrac{1}{\partial_J H}\right|_{\varepsilon=0}=\dfrac{-\partial_J H_5(I,J_0,\vec{\varphi})}{\omega_2^2}.
\end{align}
When $l\geq4$, using the Faà di Bruno's formula again, we obtain
\begin{align}
	\dfrac{d^l}{d\varepsilon^l}\dfrac{1}{\partial_{J} H}=&\dfrac{-1}{(\partial_{J}H)^2}\dfrac{d^l}{d\varepsilon}\partial_{J}H+\sum_{\gamma\in\widetilde{\Gamma}_l}(-1)^{|\gamma|}\dfrac{l!|\gamma|!}{\gamma!(\partial_{J}H)^{|\gamma|+1}}\prod_{k=1}^{l-3}\left(\dfrac{\tfrac{d^{k+1}}{d\varepsilon}\partial_{J}H}{(k+1)!}\right)^{\gamma_k}.
\end{align}
Summarizing,
\begin{align}
	\dfrac{1}{\partial_J H(I,J,\vec{\varphi})}=&\dfrac{1}{\omega_2}-\varepsilon^2\dfrac{\partial_J H_1(I,J_0)}{\omega_2^2}-\varepsilon^3\dfrac{\partial_J H_5(I,J_0,\vec{\varphi})}{\omega_2^2}\\&\ +\sum_{l=4}^{N}\dfrac{\varepsilon^l}{l!}\left(\left.\dfrac{d^l}{d\varepsilon}\dfrac{1}{\partial_{J} H}\right|_{\varepsilon=0}\right)+\Theta(\varepsilon^{N+1}),
\end{align}
where $\left.\dfrac{d^l}{d\varepsilon^l}\dfrac{1}{\partial_{J} H}\right|_{\varepsilon=0}$ has the following expression for $l\geq4$:
\begin{align}
	\left[\dfrac{-l!}{\omega_2^2}\partial_JH_{l+2}+\sum_{\gamma\in\widetilde{\Gamma}_l}(-1)^{|\gamma|}\dfrac{l!|\gamma|!}{\gamma!}\dfrac{1}{\omega_2^{|\gamma|+1}}\prod_{k=1}^{l-3}\left(\dfrac{\partial_{J}H_{k+2}}{k!}\right)^{\gamma_k}\right](I,J_0,\vec{\varphi}).
\end{align}
Finally, let $h_3$ be defined as in the beginning of the proof and set
\begin{equation*}
	h_l:=(l!)^{-1}\left.\frac{d^l}{d\varepsilon^l}\frac{1}{\partial_{J} H}\right|_{\varepsilon=0},\quad \hbox{for }l\geq4,
\end{equation*}
which completes the proof. 
\cqfd

To conclude, we now prove Proposition~\ref{Pro-Hyp0} by combining these auxiliary results. Indeed, from Propositions~\ref{Pro-Hl-numerator} and \ref{Pro-Hl-denominator} we can write
\begin{align*}
	\sum\limits_{l=5}^{N}\varepsilon^{l-2}\dfrac{\partial_\varphi H_l(I,J,\vec{\varphi})}{\partial_J H(I,J,\vec{\varphi})}&=\sum\limits_{l=3}^{N}\varepsilon^l h^\varphi_l\left(\dfrac{1}{\omega_2}-\varepsilon^2\dfrac{\partial_J H_1}{\omega_2^2}+\sum_{l=3}^{N}\varepsilon^lh_l\right)+\Theta(\varepsilon^{N+1}),\\
	&=\sum\limits_{l=3}^{N}\varepsilon^l F_l+\Theta(\varepsilon^{N+1}),
\end{align*}
from which it easily follows that
\begin{align*}
	F_l=-\dfrac{\hh_l^\varphi}{\omega_2}+\text{l.o.t.}\ \hbox{ for }l\in\{3,4,5\},
\end{align*}
and,
\begin{align*}
	F_l=\left(\dfrac{-1}{\omega_2}\right)^{\left\lfloor\tfrac{l}{3}\right\rfloor}\left[\hh_l^\varphi+\sum_{k\geq3}^{l-3}\hh_{l-k}^\varphi\hh_k\right]+\text{l.o.t.}\ \hbox{ for }l\geq6.
\end{align*}
On the one hand, by Proposition~\ref{Pro-Hl-numerator}, 
\begin{equation}
\hh^\varphi_l\in\mathcal{L}_{[l]_{\ZZ_3},\varphi}^{5+[l]_{\ZZ_3}}\left[\left\lfloor\tfrac{l}{3}\right\rfloor,\left\lfloor\tfrac{l}{3}\right\rfloor-1\right],\quad\hbox{for all }l\geq3.
\end{equation}
On the other hand, concerning the order of $\hh_{l-k}^\varphi\hh_k$, we have the following:
\begin{align*}
	\dd[\hh_{l-k}^\varphi\hh_k]&=l+2\left(\left\lfloor\frac{l-k}{3}\right\rfloor+\left\lfloor\frac{k}{3}\right\rfloor\right)=l+2\left(\left\lfloor\frac{l}{3}\right\rfloor+\left\lfloor\frac{[l]_{\ZZ_3}-[k]_{\ZZ_3}}{3}\right\rfloor\right)\leq\dd[\hh_l^\varphi],
\end{align*}
for all $3\leq k\leq l-3$, and $l\geq6$. In particular, 
\begin{equation*}
    \dd[\hh_{l}^\varphi]=\dd[\hh_{l-k}^\varphi\hh_k]\Leftrightarrow [l]_{\ZZ_3}\geq[k]_{\ZZ_3}.
\end{equation*}
Consequently, the function $\ff_l$ from Proposition~\ref{Pro-Hyp0} is characterized as follows: 
\begin{equation*}
	\ff_l=\hh_l^\varphi,\ \ \hbox{for }l\in\{3,4,5\},\quad\hbox{and}\quad\ff_l=\hh_l^\varphi+\sum_{j=0}^{[l]_{\ZZ_3}}\sum_{k=1}^{\left\lfloor\tfrac{l}{3}\right\rfloor-1}\hh_{l-j-3k}^\varphi\hh_{j+3k},\hbox{ for }l\geq6,
\end{equation*}
which immediately yields
\begin{equation}    \ff_l\in\mathcal{L}_{[l]_{\ZZ_3},\varphi}^{5+[l]_{\ZZ_3}}\left[\left\lfloor\tfrac{l}{3}\right\rfloor,\left\lfloor\tfrac{l}{3}\right\rfloor-1\right],\quad\hbox{for all }l\geq3.
\end{equation}
Since the procedure is analogous for the function $\g_l$, this completes the proof.

\section{Proof of Proposition~\ref{Pro-Hyp1}}\label{Section-Proof2}
Using Proposition~\ref{Pro-Hyp0} we can write
\begin{equation*}
	I(\theta)=I_0+\sum_{l\geq3}\varepsilon^l\int_{0}^{\theta} F_l(\theta')d\theta',
\end{equation*}
and 
\begin{equation*}
	\varphi(\theta)=\varphi_0+nm^{-1}\theta+\varepsilon^2\int_{0}^{\theta}\Lambda(I(\theta'))d\theta'+\sum_{l\geq3}\varepsilon^l\int_{0}^{\theta} G_l(\theta')d\theta',
\end{equation*}
where, to simplify notation, we denote by $\theta'$ the argument $(I(\theta'),J_0(I(\theta')),\varphi(\theta'),\theta')$ in every function of the set $\{F_l,\ G_l,\ l\geq3\}$. In addition, we also write the total derivative operators with respect to $\varepsilon$ as
\begin{equation*}
	\deps^k:=\left.\dfrac{d^k}{d\varepsilon^k}\ \right|_{\varepsilon=0},\ \hbox{for all }k\geq1,
\end{equation*}
while $\deps^0$ denotes the identity operator.

It is immediate that $\deps I=\deps^2 I\equiv0$. Hence, using the generalized Leibniz rule we obtain
\begin{align*}
	\deps^k I=\sum_{l\geq3}\sum_{j=0}^{k}\begin{pmatrix}
		k\\j
	\end{pmatrix}\left(\deps^j\varepsilon^l\right)\int_{0}^{\theta}\deps^{k-j}F_l(\theta')d\theta'=\sum_{l=3}^{k}\dfrac{k!}{(k-l)!}\int_{0}^{\theta}\deps^{k-l}F_l(\theta')d\theta',
\end{align*}
for all $k\geq3$. On the other hand,
\begin{equation*}
	\deps^k\int_{0}^{\theta}\Lambda(I(\theta'))d\theta'=\Lambda_V\int_{0}^{\theta}\deps^{k}Id\theta',
\end{equation*}
which leads to
\begin{equation*}
	\varepsilon^2\int_{0}^{\theta}\Lambda(I(\theta'))d\theta'=\varepsilon^2\Lambda(I_0)\theta+\sum_{k\geq3}^{N}\varepsilon^k\Lambda_V\int_{0}^{\theta}\deps^{k-2}Id\theta'+\Theta(\varepsilon^{N+1}).
\end{equation*}
Therefore, similarly to $I(\theta)$, we have:
\begin{align*}
	\deps^k \varphi=\Lambda_V\int_{0}^{\theta}\deps^{k-2}Id\theta'+\sum_{l=3}^{k}\dfrac{k!}{(k-l)!}\int_{0}^{\theta}\deps^{k-l}G_l(\theta')d\theta',
\end{align*}
for all $k\geq3$. In order to compute the derivatives $\deps^{k-l}F_l$ and $\deps^{k-l}G_l$, since both $F_l$ and $G_l$ satisfy \eqref{J-Constants}, we can write them as follows:
\begin{equation*}
	\sum\limits_{\gamma\in\Gamma_{\hat{l}}}f_\gamma(\varphi)\mathbf{C}[\gamma,I] :=\sum\limits_{\gamma\in\Gamma_{\hat{l}}} i^{\gamma_3+\gamma_4}e^{-i\varphi(\gamma_1-\gamma_3)}e^{-i\theta(\gamma_2-\gamma_4)}\mathbf{C}[\gamma,I].
\end{equation*}
Here we are identifying $\hat{l}:=\dd[F_l]$, where $\dd[\cdot]$ denotes the polynomial degree introduced before Lemma~\ref{J-Lemma}, while the structure of the constants $\mathbf{C}[\gamma,I]$ is given in \eqref{J-Constants}. Moreover, since the procedure is analogous in both cases, we present the computations only for $F_l$.

First, by applying the generalized Leibniz rule, we obtain:
\begin{align}\label{der-Fl}
	\deps^kF_l(\theta)=&\sum\limits_{\gamma\in\Gamma_{\hat{l}}}\sum_{j=0}^{k}\begin{pmatrix}
		k\\j
	\end{pmatrix}\deps^jf_\gamma(\varphi)\deps^{k-j}\mathbf{C}[\gamma,I].
\end{align}
On the other hand, observe that
\begin{equation*}
	\deps I=\deps^2 I=\deps \varphi\equiv0,\quad\hbox{and}\quad\deps^2\varphi=\Lambda(I_0)\theta,
\end{equation*}for all $\gamma\in\Gamma_{\hat{l}}$, with $l\geq3$, and all $I_0\in(0,d_I)$. Then, it follows that 
\begin{align}\label{FirstDerivativesEpsilon}
	\deps f_\gamma(\varphi)=\deps \mathbf{C}[I,\gamma]=\deps^2 \mathbf{C}[I,\gamma]=0,\ \hbox{and }\deps^2 f_\gamma(\varphi)=\partial_\varphi f_\gamma(\varphi^0)\Lambda(I_0)\theta,
\end{align}
where $\varphi^0=\varphi_0+nm^{-1}\theta$. Consequently, $\deps F_l=\deps G_l\equiv0$, and
\begin{equation}\label{SecondOrderDerivativeFl}
	\deps^2F_l(\theta)=\Lambda(I_0)\theta\sum\limits_{\gamma\in\Gamma_{\hat{l}}}\partial_\varphi f_\gamma(\varphi^0)\mathbf{C}[\gamma,I_0],\ \hbox{for all }l\geq3.
\end{equation}
In addition, 
\begin{equation}\label{der3-fgammaC}
	\deps^3f_\gamma(\varphi)=\partial_\varphi f_\gamma(\varphi^0)\deps^3\varphi,\quad \deps^{3}\mathbf{C}[\gamma,I]=\partial_I\mathbf{C}[\gamma,I_0]\deps^{3}I,
\end{equation}while, when $j\geq4$, the Faà di Bruno's formula gives the following:
\begin{align}\label{der-fgamma}
	\deps^jf_\gamma(\varphi)&=\partial_\varphi f_\gamma(\varphi^0)\deps^j\varphi+\sum_{\widetilde{\gamma}\in\widetilde{\Gamma}_{j}}\dfrac{j!}{\widetilde{\gamma}!}\partial_\varphi^{|\widetilde{\gamma}|}f_\gamma(\varphi^0)\prod_{i=1}^{j-3}	\left(\dfrac{\deps^{i+1}\varphi}{(i+1!)}\right)^{\widetilde{\gamma}_i},
\end{align}
and
\begin{align}\label{der-C}
	\deps^{j}\mathbf{C}[\gamma,I]&=\partial_I\mathbf{C}[\gamma,I_0]\deps^{j}I+\sum_{\widetilde{\gamma}\in\widetilde{\Gamma}_{j}}\dfrac{j!}{\widetilde{\gamma}!}\partial_I^{|\widetilde{\gamma}|}\mathbf{C}[\gamma,I_0]\prod_{i=1}^{j-3}	\left(\dfrac{\deps^{i+1}I}{(i+1!)}\right)^{\widetilde{\gamma}_i},
\end{align}
where the notation $\widetilde{\Gamma}_j$ is introduced in~\eqref{SetGamma}. 

Finally, using these expressions we establish the following lemma, which, when combined with Corollary~\ref{Cor-Hyp} and taking $\theta=2\pi m$, directly implies Proposition~\ref{Pro-Hyp1}.
\begin{Lem}
	Given $V\in\mathcal{V}_{N_0}$, consider a level set \eqref{levelSet} of sufficiently small energy $h>0$, and fix an arbitrary point $(I_0,\varphi_0)\in(0,d_I)\times\TT$. Then, for all $k\geq3$, the derivatives $\deps^k I$ and $\deps^k \varphi$ are polynomials of degree $k+2\left\lfloor\tfrac{k}{3}\right\rfloor$ in the angular variables defined in \eqref{levelSet-AngularVariables}. \\
	More concretely, when $k\in\{3,4,5\},$ we have
	\begin{align*}
		\deps^k I=k!\int_0^{\theta}F_k\left(\hat{X}(\theta')\right)d\theta'+\hbox{ l.o.t.}\quad\hbox{and}\quad
		\deps^k \varphi=k!\int_0^{\theta}G_k\left(\hat{X}(\theta')\right)d\theta'+\hbox{ l.o.t.}
	\end{align*}
	with $\hat{X}(\theta)=(I_0,J_0(I_0),\varphi_0+nm^{-1}\theta,\theta)$. Moreover, when $k\geq6$, the expressions take the form:
	\begin{align*}
		\deps^k I=k!\int_0^{\theta}F_k\left(\hat{X}(\theta')\right)d\theta'+\sum_{j=0}^{[k]_{\ZZ_3}}\sum_{l=1}^{\left\lfloor\tfrac{k}{3}\right\rfloor-1}k!\int_{0}^{\theta} \dfrac{ \deps^{k-3l-j} F_{3l+j}\left(\hat{X}(\theta')\right)}{(k-3l-j)!}d\theta'+\hbox{ l.o.t.}
	\end{align*}
	and
	\begin{align*}
		\deps^k \varphi=k!\int_0^{\theta}G_k\left(\hat{X}(\theta')\right)d\theta'+\sum_{j=0}^{[k]_{\ZZ_3}}\sum_{l=1}^{\left\lfloor\tfrac{k}{3}\right\rfloor-1}k!\int_{0}^{\theta}\dfrac{\deps^{k-3l-j} G_{3l+j}\left(\hat{X}(\theta')\right)}{(k-3l-j)!} d\theta'+\hbox{ l.o.t.}
	\end{align*}
\end{Lem}
\proof
As usual, denoting the Kronecker's delta by $\delta_{ij}$, when $k\in\{3,4,5\}$, we have 
\begin{align*}
	\deps^kI=k!\int_{0}^{\theta} F_k\left(\hat{X}(\theta')\right)d\theta'+\delta_{k5}\dfrac{k!}{2!}\int_{0}^{\theta} \deps^2 F_{k-2}\left(\hat{X}(\theta')\right)d\theta',
\end{align*} and
\begin{align*}
	\deps^k\varphi=k!\int_{0}^{\theta} G_k\left(\hat{X}(\theta')\right)d\theta'+\delta_{k5}\left(\Lambda_V\int_{0}^{\theta}\deps^3 Id\theta'+\dfrac{k!}{2!}\int_{0}^{\theta} \deps^2 G_{k-2}\left(\hat{X}(\theta')\right)d\theta'\right). 
\end{align*}
Moreover, for all $\gamma\in\Gamma_{\hat{l}}$ with $l\geq3$, it is satisfied that
\begin{align*}
	\dd\left[\partial_\varphi f_\gamma\Cc[\gamma,I]\right]=\dd\left[f_\gamma\Cc[\gamma,I]\right]=\dd[F_l]=\dd[G_l],
\end{align*}
which, combined with \eqref{SecondOrderDerivativeFl}, leads to
\begin{align}\label{Fl-equality}
	\dd\left[\deps^2 F_l\right]=\dd\left[\deps^2 G_l\right]=\dd\left[F_l\right]=\dd\left[G_l\right],\quad \hbox{for all }l\geq3.
\end{align}
Note that the degree of a polynomial with respect to the angular variables in \eqref{levelSet-AngularVariables} remains invariant under integration in $\theta$, which implies the lemma for these values of $k$. Complementary, when $k\geq6$, we have the expressions:
\begin{align*}
	&\deps^k I=k!\int_{0}^{\theta} F_k\left(\hat{X}(\theta')\right)d\theta'+\sum_{l=3}^{k-2}\dfrac{k!}{(k-l)!}\int_{0}^{\theta}\deps^{k-l}F_l\left(\hat{X}(\theta')\right)d\theta',\\ &\deps^k\varphi=k!\int_{0}^{\theta} G_k\left(\hat{X}(\theta')\right)d\theta'+\Lambda_V\int_{0}^{\theta}\deps^{k-2} Id\theta'+\sum_{l=3}^{k-2}\dfrac{k!}{(k-l)!}\int_{0}^{\theta}\deps^{k-l}G_l\left(\hat{X}(\theta')\right)d\theta'.
\end{align*}
Then, since $F_l$ and $G_l$ share the same structure (and therefore do so $\deps^k G_l$ and $\deps^k F_l$, for all $l,k$), observe that the proof ends if the following holds, for all $k\geq6$, and $3\leq l\leq k-2$:
\begin{align}\label{LastInequality}
	\left\{\begin{array}{l}
		\dd\left[F_k\right]=\dd\left[\deps^{k-l}F_{l}\right]\ \hbox{when }[l]_{\ZZ_3}\leq[k]_{\ZZ_3},\vspace{1mm}\\ 
		\dd\left[F_k\right]>\dd\left[\deps^{k-l}F_{l}\right]\ \hbox{otherwise.}
	\end{array}\right.
\end{align}
The case $l=k-2$ is already proved by \eqref{Fl-equality}. Additionally, let us emphasize that
\begin{align*}
	\dd\left[\deps f_\gamma\right]=\dd\left[\deps \mathbf{C}[\gamma,I]\right]=\dd\left[\deps^2\mathbf{C}[\gamma,I]\right]=0,\ \ \dd\left[\deps^2 f_\gamma\right]=\dd\left[f_\gamma\right]=\dd\left[F_l\right],
\end{align*}
for all $\gamma\in\Gamma_l$, with $l\geq3$. Accordingly, from \eqref{der-Fl} we can write the following:
\begin{align}\label{der-Fl2}
	\deps^{k-l}F_l(\theta)=&\sum\limits_{\gamma\in\Gamma_{l}}\left[\mathbf{C}[\gamma,I_0]\deps^{k-l}f_\gamma(\varphi)+f_\gamma(\varphi)\deps^{k-l}\mathbf{C}[\gamma,I]\right]\\&\ +\delta_{5,k-l}\tfrac{5!}{2}\sum\limits_{\gamma\in\Gamma_{l}}\deps^2f_\gamma(\varphi)\deps^{3}\mathbf{C}[\gamma,I],\quad\hbox{when }3\leq k-l\leq5.
\end{align}
Then, using \eqref{der3-fgammaC} we obtain that
\begin{align*}
	\dd\left[\mathbf{C}[\gamma,I_0]\deps^{3}f_\gamma(\varphi)+f_\gamma(\varphi)\deps^{3}\mathbf{C}[\gamma,I]\right]=\dd\left[\deps^{3}f_\gamma(\varphi)\right]=\dd\left[f_\gamma(\varphi)\deps^{3}\mathbf{C}[\gamma,I]\right],
\end{align*}
for any $\gamma\in\Gamma_3$, and thus
\begin{equation*}
	\dd[\deps^6 I]=\dd[\deps^6\varphi]=\dd[\deps^3 F_3]=\dd[\deps^3 G_3]=10.
\end{equation*}
In particular, this proves the lemma for the case $k=6$. \\
On the other hand, when $(k-l)\geq6$:
\begin{align}\label{der-Fl3}
	\deps^{k-l}F_l(\theta)=&\sum\limits_{\gamma\in\Gamma_{l}}\left[\mathbf{C}[\gamma,I_0]\deps^{k-l}f_\gamma(\varphi)+f_\gamma(\varphi)\deps^{k-l}\mathbf{C}[\gamma,I]\right]\\&\ +\frac{(k-l)!}{2}\sum\limits_{\gamma\in\Gamma_{l}}\deps^2f_\gamma(\varphi)\deps^{k-l-2}\mathbf{C}[\gamma,I]\\&+\sum\limits_{\gamma\in\Gamma_{l}}\left[\sum_{j=3}^{k-l-3}\begin{pmatrix}
		k-l\\ j
	\end{pmatrix}\deps^j f_\gamma \deps^{k-l-j}\mathbf{C}[\gamma,I]\right].
\end{align}
Therefore, it is necessary to delve into the structure of $\deps^j f_\gamma$ and $\deps^j \Cc[\gamma,I]$ to prove \eqref{LastInequality}. To this end, let us assume that the statements are true up to some fixed $\kk\geq6$, i.e.,
\begin{align*}
	\dd[\deps^lI]=\dd[\deps^l\varphi]=l+2\left\lfloor\tfrac{l}{3}\right\rfloor,\quad\hbox{for all }3\leq l\leq \kk.
\end{align*}
As a consequence, the set $\widetilde{\Gamma}^M_{l}$ defined by \eqref{J-MaxSet} coincides for $J_l$, $\deps^lI$ and $\deps^l\varphi$, that is: for any $4\leq l\leq \kk$, the set from \eqref{J-MaxSet} can be equivalently defined as
\begin{align*}
	\widetilde{\Gamma}^M_{l}&=\left\{\gamma\in\widetilde{\Gamma}_{l}:{\dd}\left[\prod_{j}(\deps^{j+1}I)^{\gamma_j}\right]={\dd}[\deps^lI]\right\}\\&=\left\{\gamma\in\widetilde{\Gamma}_{l}:{\dd}\left[\prod_{j}(\deps^{j+1}\varphi)^{\gamma_j}\right]={\dd}[\deps^l\varphi]\right\}.
\end{align*}
Then, Corollary~\ref{J-Corollary} can be applied in each case, which implies that $\widetilde{\Gamma}^M_{4}=\widetilde{\Gamma}^M_{5}=\emptyset$. Moreover, for all $\gamma\in\Gamma_{\hat{l}}$, with $l\geq3$, and all $4\leq j\leq \kk$, from \eqref{der-fgamma} and \eqref{der-C} we can write
\begin{align*}
	\deps^jf_\gamma(\varphi)&=\partial_\varphi f_\gamma(\varphi^0)\deps^j\varphi+\sum_{\widetilde{\gamma}\in\widetilde{\Gamma}^M_{j}}\dfrac{j!}{\widetilde{\gamma}!}\partial_\varphi^{|\widetilde{\gamma}|}f_\gamma(\varphi^0)\prod_{i=1}^{j-3}	\left(\dfrac{\deps^{i+1}\varphi}{(i+1!)}\right)^{\widetilde{\gamma}_i}+\hbox{ l.o.t.}
\end{align*}
and
\begin{align*}
	\deps^{j}\mathbf{C}[\gamma,I]&=\partial_I\mathbf{C}[\gamma,I_0]\deps^{j}I+\sum_{\widetilde{\gamma}\in\widetilde{\Gamma}^M_{j}}\dfrac{j!}{\widetilde{\gamma}!}\partial_I^{|\widetilde{\gamma}|}\mathbf{C}[\gamma,I_0]\prod_{i=1}^{j-3}	\left(\dfrac{\deps^{i+1}I}{(i+1!)}\right)^{\widetilde{\gamma}_i}+\hbox{ l.o.t.}
\end{align*}
Hence,
\begin{align*}
	\dd\left[\deps^j f_\gamma\right]=\dd\left[f_\gamma\right]+\dd\left[\deps^j \varphi\right],\ \ \dd\left[\deps^j \mathbf{C}[\gamma,I]\right]=\dd\left[\deps^jI\right],
\end{align*}
which leads to
\begin{align}\label{derivatives-fgammaC}
	\dd\left[\deps^j f_\gamma\right]=l+j+2\left(\left\lfloor\tfrac{l}{3}\right\rfloor+\left\lfloor\tfrac{j}{3}\right\rfloor\right),\quad \dd\left[\deps^j \mathbf{C}[\gamma,I]\right]=j+2\left\lfloor\tfrac{j}{3}\right\rfloor,
\end{align}
for all $\gamma\in\Gamma_{\hat{l}}$, with $l\geq3$, and all $3\leq j\leq \kk$. Using this, for any $6\leq k\leq\kk$, it is immediate that $$\dd\left[\deps^{k-l}f_\gamma(\varphi)\right]=\dd\left[f_\gamma(\varphi)\deps^{k-l}\mathbf{C}[\gamma,I]\right],\hbox{ for all }\gamma\in\Gamma_{\hat{l}},\hbox{ with }3\leq l\leq k-3. $$ Moreover,
\begin{align*}
	\left\{\begin{array}{l}
		\dd\left[\deps^{k-l}f_\gamma(\varphi)\right]=\dd[\deps^k I]\ \hbox{ when }[l]_{\ZZ_3}\leq[k]_{\ZZ_3},\vspace{1mm}\\
		\dd\left[\deps^{k-l}f_\gamma(\varphi)\right]<\dd[\deps^k I]\ \hbox{ otherwise,}
	\end{array}\right.
\end{align*}
and
\begin{align*}
	\left\{\begin{array}{l}
		\dd\left[\deps^{2}f_\gamma(\varphi)\deps^{k-l-2}\mathbf{C}[\gamma,I]\right]=\dd[\deps^k I]\ \hbox{ when }k-l=5,\vspace{1mm}\\
		\dd\left[\deps^{2}f_\gamma(\varphi)\deps^{k-l-2}\mathbf{C}[\gamma,I]\right]<\dd[\deps^k I]\ \hbox{ otherwise,}
	\end{array}\right.
\end{align*}
for all $\gamma\in\Gamma_{\hat{l}}$ with $3\leq l\leq k-3$. As a consequence, by considering the function $\deps^{k-l}F_l$ defined in \eqref{der-Fl2}, it follows that \eqref{LastInequality} is satisfied for all $3\leq k-l\leq 5$, thus proving the lemma for the cases $k\in\{7,8\}$. Furthermore, the case $k=9$ is also proved by taking $j=3$ in \eqref{der-Fl3}, leading to
\begin{align*}
	\dd\left[\deps^{3}f_\gamma(\varphi)\deps^{3}\mathbf{C}[\gamma,I]\right]=\dd\left[\deps^9 I]\right]=15,\quad\hbox{for all }\gamma\in\Gamma_{\hat{l}},\ \hbox{ with }l=3.
\end{align*}
Therefore, we can consider $9\leq k\leq\kk$. In that case, let $3\leq l\leq k-6$ and $3\leq j\leq k-l-3$. Then,
\begin{align*}
	\dd\left[\deps^j f_\gamma \deps^{k-l-j}\mathbf{C}[\gamma,I]\right]\leq k+2\left(\left\lfloor\tfrac{k}{3}\right\rfloor+\left\lfloor\tfrac{l}{3}\right\rfloor+\left\lfloor\tfrac{j}{3}\right\rfloor-\left\lfloor\tfrac{l+j}{3}\right\rfloor\right)\leq k+2\left\lfloor\tfrac{k}{3}\right\rfloor,
\end{align*}
where the first inequality becomes an equality when $[l+j]_{\ZZ_3}\leq[k]_{\ZZ_3}$ and the second when $[l]_{\ZZ_3}=[j]_{\ZZ_3}$. By applying this in \eqref{der-Fl3}, together with the above estimates, it is straightforward to prove \eqref{LastInequality} for all $6\leq k\leq\kk$ and $3\leq l\leq k-2$. 

To conclude, take $k=\kk+1$ and consider the function $d_\varepsilon^{k-l}F_l$ defined in \eqref{der-Fl2} and \eqref{der-Fl3}. Observe that every derivative $\deps^j$ appearing in the right-hand sides of these expressions satisfies that $j\leq \kk-2$. In particular, all these terms satisfies the identities from \eqref{derivatives-fgammaC}. Therefore, starting from \eqref{derivatives-fgammaC}, the analogous procedure yields the corresponding bounds for the case $\kk+1$, thus proving \eqref{LastInequality} and, consequently, the lemma.
\cqfd

\section*{Acknowledgments} This work has received funding from the European Research Council (ERC) under the European Union's Horizon 2020 research and innovation programme through the grant agreement~862342 (A.E.). The authors are partially supported by the grants CEX2023-001347-S, RED2022-134301-T and PID2022-136795NB-I00 (A.E., and D.P.-S.) funded by MCIN/AEI/10.13039/501100011033, and Ayudas Fundaci\'on BBVA a Proyectos de Investigaci\'on Cient\'ifica 2021 (D.P.-S.).  M.G. is partially supported by the Spanish MCIN/AEI grant JDC2022-049490-I, project PID2021-128418NA-I00.

\appendix
\section{The class $\mathcal{V}_{N_0}$ is dense in $ C^\infty(\TT^2)$}\label{Appendix-Denseclass}
For any fixed integer $N_0\geq2$, we begin by recalling the definition \eqref{coprimes}: a potential $V\in\mathcal{V}_{0}$ satisfying  \eqref{nonresonant}-\eqref{Isochronic} belongs to $\mathcal{V}_{N_0}$ only if $\tfrac{\omega_1}{\omega_2}=\tfrac{n}{m}\in\mathbb{Q}$, where
\begin{equation}\label{coprimes-2}
	n,m\hbox{ are coprime, }n\geq2,\ m\geq N_0,\hbox{ and }[n+m]_{\ZZ_5}\leq2.
\end{equation}

It is well known that the set of rational numbers with both the numerator and the denominator prime is dense in $(0,+\infty)$ (see, for instance, \cite[pg. 165]{Sierpinski}). Then, the claim is proved if any positive rational number with prime coefficients can be approximated by a sequence $\{\tfrac{n_k}{m_k}\}\subset\mathbb{Q}$ satisfying \eqref{coprimes-2}. To this end, fix two prime numbers $n,m$, and define the sequence
\begin{equation*}
n_k=5kn+1,\ m_k=5km,\ \ \hbox{for every prime number $k\geq \tfrac{N_0}{5m}$}.
\end{equation*}
Since $m_k\geq N_0$ for all $k$, we only have to check that they are coprime, which is equivalent to see that $\{5,k,m\}$ are not divisors of the integer $n_k$.

First, it is clear that $n_k$ is not divisible by $5$ or $n$, for any $k$ prime. However, assume that there exists a sequence of prime numbers $\{k_j\}\nearrow+\infty$ such that $n_{k_j}$ is divisible by $m$, for all $j\in\mathbb{N}$. Then, consider the sequence
\begin{equation*}
	\tilde{n}_{k_j}=5k_jn+2,\ \tilde{m}_{k_j}=m_{k_j}.
\end{equation*}
Since they are coprime, the proof is complete.

\end{document}